\documentclass[a4paper,10pt]{article}

\usepackage{geometry}
\usepackage{amsmath}
\usepackage{amsfonts}
\usepackage{amsthm}

\usepackage[latin1]{inputenc}  
\usepackage[T1]{fontenc}       
\usepackage{amsfonts}
\usepackage{amsthm}
\usepackage{amssymb}
\usepackage[pdftex]{graphicx}
\usepackage{dsfont}

\usepackage{hyperref} 
\hypersetup{ 
colorlinks=true, 
breaklinks=true, 
urlcolor= blue, 
linkcolor= red, 
} 

\usepackage{psfrag}

\newtheorem{theo}{Theorem}[section]
\newtheorem{def1}[theo]{Definition}
\newtheorem{cor}[theo]{Corollary}
\newtheorem{prop}[theo]{Proposition}
\newtheorem*{rem}{Remark}
\newtheorem{lem}[theo]{Lemma}

\def\eps{{\varepsilon}}
\def\a{{\alpha}}

\def\R{{\mathbb R}}
\def\Sp{{\mathbb S}}

\def\Sc{{\mathcal S}}
\def\N{{\mathbb N}}

\def\Ot{{\mathbb O}}
\def\1{{\mathds 1}}

\title{Radiation condition at infinity for the high-frequency Helmholtz equation: optimality of a non-refocusing criterion}

\author{ \footnote{IRMAR, UMR 6625-CNRS, Campus de Beaulieu, 35042 Rennes, France, E.mail: francois.castella@univ-rennes1.fr} Fran\c{c}ois Castella and \footnote{LATP, UMR 7353-CNRS, CMI Technop\^{o}le Ch\^ateau-Gombert, 39 sur F. Joliot Curie, 13453 Marseille Cedex 13, FRANCE, France, E.mail: aurelien.klak@cmi.univ-mrs.fr} Aur\'elien Klak}

\date{\today}

\begin{document}
\maketitle 


\begin{abstract}
We consider the high frequency Helmholtz equation with a variable refraction index $n^2(x)$ ($x \in \R^d$),
supplemented with a given high frequency source term supported near the origin $x=0$.
A small absorption parameter $\alpha_{\varepsilon}>0$ is added, which
somehow prescribes a radiation condition at infinity for the considered Helmholtz equation.
The semi-classical parameter is $\varepsilon>0$. We let
$\eps$ and $\a_\eps$ go to zero
{\em simultaneaously}. We study the question whether the indirectly prescribed radiation condition at infinity
is satisfied {\em uniformly} along the asymptotic process $\eps \to 0$, or, in other words,
whether the conveniently rescaled solution to the considered equation goes to the {\em outgoing} solution
to the natural limiting Helmholtz equation.

This question 
has been previously studied by the first autor in \cite{MR2139886}. In  \cite{MR2139886},
it is proved that the radiation condition is indeed satisfied uniformly in $\eps$, 
provided the refraction index satisfies a specific {\em non-refocusing condition}, a condition that is first pointed out in this reference.
The non-refocusing condition requires, in essence, that the rays of geometric optics naturally associated with the high-frequency Helmholtz operator, and that are sent from the origin $x=0$ at time $t=0$, should not refocus at some later time $t>0$ near the origin again. 

In the present text we show the {\em optimality}
of the above mentionned non-refocusing condition, in the following sense. We exhibit  a refraction index
which {\em does} refocus the rays of geometric optics sent from the origin near the origin again, and, on the other hand, we completely compute the asymptotic behaviour of the solution to the associated  Helmholtz equation: we show that
the limiting solution {\em does not} satisfy the natural radiation condition at infinity. More precisely, we
show that the limiting solution is a {\em perturbation} of the outgoing solution to the natural limiting Helmholtz equation, and
that the perturbing term explicitly involves the contribution of the rays radiated from the origin which go back to the origin.
This term is also conveniently modulated by a phase factor, which turns out to be the action along the  above rays of 
the hamiltonian associated with the semiclassical Helmholtz equation.
\end{abstract}



\tableofcontents


\section{Introduction}



\subsection{General introduction}


In this article, we study the convergence as $\varepsilon$ approaches $0$ of $w^{\varepsilon}$, solution  to the following rescaled Helmholtz equation
\begin{equation}
\label{eqHelm}
i\varepsilon\alpha_{\varepsilon} \, w_\varepsilon(x)+\frac{\Delta_{x}}{2}w_\varepsilon(x)+n^2(\varepsilon x)w_\varepsilon(x)=S(x),\quad x\in\R^d \quad (d\geq3).
\end{equation}
Here $\alpha_{\varepsilon}$ is an absorption parameter, $n^2(x)$ is a space-dependent refraction index\footnote{Here and below we use the standard notation $n^2(x)$, a squared term, assuming in doing so that the corresponding term is everywhere non-negative. This is a harmless abuse of notation, since the refraction index $n^2(x)$ that is eventually chosen
in our analysis is negative for certain values of $x$. The reader may safely skip this fact, since the Helmholtz equation also arises in the spectral analysis of Schr\"odinger operators, where the refraction index becomes $E-V(x)$ where $E$ is an energy and $V(x)$ is a space-dependent potential, and the term $E-V(x)$ may change sign in that context.}
and $S(x)$ is a given and smooth source term. In the sequel, we assume the following:
\begin{itemize}
\item The absorption parameter $\alpha_{\varepsilon}$ satisfies\footnote{The limiting case $\a_\eps=0^+$ can be considered along our analysis, see below.}
\begin{equation*}
\alpha_{\varepsilon} >0, \qquad \a_\eps \mathop{\longrightarrow}_{\eps\to 0} 0.
\end{equation*}
\item The smooth refraction index $n^2(x) \in C^\infty(\R^d)$ is a possibly long-range perturbation of a 
positive constant $n_\infty^2>0$
at infinity, namely, for some $\rho>0$, we have
\begin{equation}\label{condpot}
\forall\, \alpha\in\N^d,\quad \exists C_\a, \quad \forall\, x\in\R^d, \quad
\Big|\partial_{x}^\alpha\left( n^2(x)-n^2_{\infty} \right)\Big|
\leq C_{\alpha}\langle x\rangle^{-\rho-\alpha},
\end{equation} 
where we denote as usual $\langle x \rangle :=(1+|x|^2)^{1/2}$.

\item
The source term $S(x)$ belongs to the Schwartz class\footnote{This assumption may be considerably relaxed at the price of some irrelevant technicalities.}  $\Sc(\R^d)$.

\end{itemize}

The question we raise is the following. Thanks to the absorption parameter $\a_\eps>0$ in (\ref{eqHelm}), the sequence of solutions $w_\eps$ is uniquely defined (see below for the limiting case $\a_\eps=0^+$).
On top of that, and as a consequence of specific homogeneous bounds obtained by Perthame and Vega in
\cite{MR1695559} (see \cite{MR2289541} for extensions by Jecko and the first author, as well as \cite{CJK}),
it is clear that the sequence $w_\eps$ is bounded in some
weighted $L^2$ space, uniformly in $\eps$. Hence the sequence $w_\eps$ possesses a limit (up to subsequences), say
in the distribution sense, and the limit $w=\lim w_\eps$ satisifies in the distribution sense the Helmholtz equation
\begin{equation}\label{eqhelmm0}
\frac{\Delta_{x}}{2} \, w +n^2(0) \, w=S,
\end{equation}
where the variable coefficients refraction index $n^2(\eps x)$ in (\ref{eqHelm}) has now coefficients frozen at the origin $x=0$.

Now, the difficulty is, the Helmholtz equation (\ref{eqhelmm0}) does not have a uniquely defined solution. At least two distinct
solutions exist, namely the outgoing solution, defined as
\begin{align}
\label{wwout}
w_{out}(x):&=\lim_{\delta\rightarrow 0^+}\left(i\delta+\frac{\Delta_{x}}{2}+n^2(0)\right)^{-1}S(x),
\end{align}
and the incoming solution, defined similarly as
$\displaystyle w_{in}=\lim_{\delta\rightarrow 0^+}\left(-i\delta+\frac{\Delta_{x}}{2}+n^2(0)\right)^{-1}S$.
Equivalently, the outgoing solution may be defined as the unique solution to the Helmholtz equation (\ref{eqhelmm0}) which
satisfies the so-called Sommerfeld radiation condition at infinity, namely
\begin{equation}\label{eqsommmerfeld}
\frac{x}{|x|}.\nabla_{x}w_{out}(x)+i \,\sqrt{2} \,  n(0) \, w_{out}(x)=O\left(\frac{1}{|x|^2}\right),\quad \text{ as } |x|\longrightarrow +\infty.
\end{equation}
This formulation means
that $w_{out}$ is required to oscillate like $w_{out} \sim \exp\left(-i \sqrt{2} \, n(0) |x| \right)/|x|$ as $|x| \to \infty$.
Similarly, the incoming solution satisfies  the following  radiation condition at infinity, namely 
$\displaystyle
( x / |x|) \, .\nabla_{x}w_{in}- i \, \sqrt{2} \, n(0) \, w_{in}=O\left(1/|x|^2\right)$, meaning that 
$w_{in} $ $\sim$ $  \exp\left(+i \sqrt{2} \, n(0) |x| \right)/|x|$ as $x \to \infty$.

\medskip

In that perspective, and due to the positive absorption parameter $\a_\eps>0$ in (\ref{eqHelm}),
it is natural to {\em expect} that the previously defined sequence $w_\eps$ goes to the {\em outgoing} solution $w_{out}$
to (\ref{eqhelmm0}).

This is the question we address here.

\medskip

It turns out that delicate analytical tools are needed to provide a
clean understanding of the phenomena at hand, and to establish whether $w_\eps \sim w_{out}$ as $\eps\to 0$.
The basic difficulty is a  conflict between a local and a global phenomenon. On the one hand, the obvious
fact that $w_\eps$ goes to a solution to (\ref{eqhelmm0}) is {\em local}: locally in $x$, {\em i.e.} in the distribution sense, the variable refraction index $n^2(\eps x)$ goes to the value $n^2(0)$ at the origin. On the other hand, the positive absorption parameter $\a_\eps>0$ in (\ref{eqHelm}) somehow asserts that $w_\eps$ is an {\em outgoing} solution to
$\Delta_x w_\eps/2+n^2(\eps x) \, w_\eps = S$, hence introducing the value at infinity
$\displaystyle
n_\infty=\lim_{x\to\infty} n(\eps x)=\lim_{x\to\infty} n(x)$, the solution $w_\eps$ should roughly oscillate
like $w_\eps\sim\exp\left(-i \, \sqrt{2} \, n_\infty |x| \right)/|x|$ at infinity. This is a {\em global} phenomenon.
Now, all this is to be compared with the fact that $w_{out}$ oscillates like
$w_{out}\sim\exp\left(-i \, \sqrt{2} \, n(0) |x| \right)/|x|$ at infinity.
Due to the fact that $n_\infty\neq n(0)$, 
the radiation condition at infinity satisfied by $w_\eps$ for any positive value $\eps>0$ is {\em a priori} incompatible with the
radiation condition at infinity satisfied by the expected limit $w_{out}$:
the radiation condition at infinity cannot be followed at once  uniformly in $\eps$, in any direct fashion (this is not in contradiction with 
the expected {\em local} convergence of $w_\eps$ towards $w_{out}$.)

\medskip

Before going further, let us mention that the above question stems from a series of articles
\cite{MR1924691}, 
\cite{MR1900556} 
about the {\em high-frequency} Helmholtz equation (Equation (\ref{eqHelm}) is
a {\em low-frequency} equation)
(see also \cite{MR2237164} and \cite{MR2296804} 
for similar considerations, in the case of a discontinuous refraction index,
as well as \cite{FixMe} 
and \cite{MR2718664} for the case of a variable absorption coefficient).
These two papers investigate the high-frequency behaviour, in terms of semi-classical measures, of high-frequency Helmholtz equations  of the form
\begin{equation}
\label{eqHelme}
i\varepsilon\alpha_{\varepsilon}u_\varepsilon(x)+
\frac{\varepsilon^2}{2} \, \Delta_{x} u_\varepsilon(x)+n^2( x)u_\varepsilon(x)=
\frac{1}{\varepsilon^{d/2}}S\left(\frac{x}{\varepsilon}\right)
\quad (x\in\R^d).
\end{equation}
The link between the low-frequency equation (\ref{eqHelm})
that is the purpose of this article, and the high-frequency equation
(\ref{eqHelme}) is provided by the following basic observation: the function  $w_\eps$ satisfies~ (\ref{eqHelm})
if and only if
the rescaled function
\begin{align}
\label{uepss}
u_\eps(x)=\frac{1}{\eps^{d/2}} \, w_\eps\left(\frac{x}{\eps}\right)
\end{align}
satisfies (\ref{eqHelme}).
In that picture, the main phenomenon to be described in~\eqref{eqHelme} is the possibility of {\em resonances}
between the high-frequency waves selected by the Helmholtz operator
$\varepsilon^2 \Delta_{x}/2 + n^2(x)$, and the high-frequency waves 
carried by the rescaled source term $\displaystyle \eps^{-d/2} \, S(x/\eps)$, both having the same wavelength $\eps$.
Amongst others, it is established in~\cite{MR1924691}, 
\cite{MR1900556} 
that the semiclassical measure associated with $u_\eps$ can be completely computed provided
$w_\eps$ indeed converges towards $w_{out}$, this latter requirement
being left as a conjecture in the cited papers.
This is the motivation for the question we address here. 

\medskip

In \cite{MR2139886}, the first positive convergence result $w_\eps\to w_{out}$ is established. 
This results requires, amongst others, a specific and original {\em non-refocusing condition} on the refraction index $n^2(x)$
(called "transversality condition" in the original paper).
This condition (see below for details) roughly asserts that the rays of geometric optics associated with the
the semi-classical Helmholtz operator $\eps^2\Delta_x/2+n^2(x)$ cannot focus at some positive time $t>0$ near the origin
$x=0$ when issued from the origin at time $t=0$. Later, X.-P. Wang and P. Zhang \cite{MR2184186} 
proved a similar, positive result, using a so-called virial assumption which is stronger than the above
non-refocusing condition. J.-F. Bony in \cite{MR2582438} 
establishes along quite different lines a positive result that is similar in spirit, requiring a weaker
non-refocusing condition.

\medskip

The goal of the present text is to prove in some sense the {\em optimality} of the non-refocusing condition
pointed out in \cite{MR2139886}.

We construct a refraction index $n^2(x)$ which violates the non-refocusing condition (rays of geometric
optics issued from the origin do refocus close to the origin at some later time), and, by explicitly computing the
asymptotic behaviour of $w_\eps$ thanks to an appropriate amplitude/phase representation developped in \cite{MR2139886},
we prove that
$$
w_\eps \mathop{\sim}_{\eps \to 0} w_{out} + \underbrace{\text{perturbation}}_{\neq 0},
$$
where the perturbation is computed as well. It explicitly involves the contribution of the rays issued from the origin
which go back to the origin at some positive time, modulated by a phase factor that is the action,
along these rays, of the hamiltonian associated with the high-frequency Helmholtz operator.


\subsection{The non-refocusing condition}


As already mentionned, the asymptotic behaviour of $w_\eps$ is dictated by that of the rescaled function
$u_\eps(x)=\eps^{-d/2} \, w_\eps(x/\eps)$. The function $u_\eps$
is $w_\eps$ rescaled at the semi-classical scale, see~\eqref{eqHelme} and \eqref{uepss}. This is translated by the following identity,
valid for any smooth test function $ \phi\in \Sc(\R^d)$, namely
\begin{equation*}
\forall\, \phi\in \Sc(\R^d), \qquad
\langle w_{\varepsilon},\phi\rangle
=
\left\langle u_{\varepsilon},\frac{1}{\varepsilon^{d/2}} \, \phi\left( \frac{x}{\varepsilon} \right) \right\rangle.\,
\end{equation*}
where we denote as usual $\langle w_{\varepsilon},\phi\rangle:=\displaystyle \int_{\R^d} w_\eps(x) \, \phi^\ast(x) \, dx$,
and $\ast$ denotes complex conjugation. In other words, the weak limit $\langle w_{\varepsilon},\phi\rangle$
of $w_\eps$ can be computed
as the weak limit {\em at the semi-classical scale} of $u^\eps$, namely the limit of
$\langle u_{\varepsilon},\varepsilon^{-d/2}\phi(x/\varepsilon)\rangle$.
This first observation is the main reason why semi-classical tools play a key role in our analysis.

Besides,
the asymptotic study of $(\ref{eqHelm})$ is done here by transforming the problem into a time-dependent problem. This approach, introduced in \cite{MR2139886},  has been used since by J.F.-Bony ($\cite{MR2582438}$) to study the Wigner measure associated to $(\ref{eqHelme})$, or by J. Royer (\cite{FixMe})
when the absorption $\alpha_{\varepsilon}$ depends on $x$. It
consists in writing the solution $w_\eps$ as the integral over the whole time of the propagator associated with $i\varepsilon\alpha_{\varepsilon}+\Delta_{x}/2+n^2(\varepsilon\, x)$, namely
\begin{equation}
\label{formintegral}
w_{\varepsilon}(x)=i\int_{0}^{+\infty}e^{-\alpha_{\varepsilon}t}e^{it\left(\frac{\Delta_{x}}{2}+n^2(\varepsilon\,x)\right)} \, S(x) \, dt.
\end{equation}
In the same way the outgoing solution can be written as
\begin{align*}
w_{out}(x):&=i\int_{0}^{+\infty}e^{it\left(\frac{\Delta_{x}}{2}+n^2(0)\right)} \, S(x) \, dt.
\end{align*}
In that picture, proving or disproving the convergence $w_\eps \sim w_{out}$ reduces to passing to the limit
in the above time integral.

Combining the two above observations, the basic first step of our analysis
consists in writing, for any given test function $\phi$, an in \cite{MR2139886},
\begin{align}
\label{basicform}
&
\nonumber
\langle w_\eps \, , \, \phi \rangle
=
\big\langle u^\eps \, , \, \eps^{-d/2} \phi(x/\eps) \big\rangle
\\
&
\qquad\quad\,\,
=
\frac{i}{\eps} \, \int_0^{+\infty}
e^{-\a_\eps \, t} \,
\Big\langle U_\eps(t) \, S_\eps \, , \, \phi_\eps\Big\rangle \, dt,
\end{align}
where we use the notation
\begin{align}
\label{resca}
&
S_\eps(x):= \frac{1}{\eps^{d/2}} \, S\left(\frac{x}{\eps} \right), \quad
\text{ and  similarly } \, 
\phi_\eps(x):= \frac{1}{\eps^{d/2}} \, \phi\left(\frac{x}{\eps} \right),
\end{align}
where  the semi-classical propagator associated with the semi-classical Hamiltonian $\eps^2 \Delta_x/2 + n^2(x)$ is
\begin{align}
\label{propa}
U_\eps(t)=\exp\left( i \, \frac{t}{\eps} \, \left( \frac{\eps^2}{2} \Delta_x + n^2(x) \right) \right)
\end{align}
It is fairly clear on formula (\ref{basicform}) that the asymptotics $\eps \to 0$ in $\langle w_\eps , \phi \rangle$ 
is dominated on the one hand by the concentration of the rescaled test function $\phi_\eps$ close
to the origin at the semi-classical scale~$\eps$, and on the other hand by the oscillations induced by
the semi-classical propagator $U^\eps(t)$ at the semi-classical scale $\eps$ as well.
The point is to measure the possible constructive intereference between both waves.

\medskip

As standard in semiclassical analysis we define the semiclassical symbol
\begin{align}
\label{ham}
h(x,\xi)=\frac{|\xi|^2}{2}-n^2(x),
\end{align}
associated with the semiclassical Schr\"odinger operator $-\frac{\varepsilon^2}{2} \, \Delta_{x}-n^2(x)$.
The semi-classical propagator $U_\eps(t)$ is known to roughly propagate
the information along the rays of geometric optics, defined as the solutions to the Hamiltonian ODE associated with $h$,
namely (see {\em e.g.} \cite{MR1735654}, \cite{MR1872698}, or \cite{MR897108})
\begin{equation}
\label{eqHJ}
\begin{cases}
\begin{array}{ll}
\displaystyle
\vspace{0.2cm}
\frac{\partial}{\partial t}X(t,x,\xi)=\Xi(t,x,\xi),&   \quad X(0,x,\xi)=x,\\
\displaystyle
\frac{\partial}{\partial t}\Xi(t,x,\xi)=\nabla_x n^2(X(t,x,\xi)),&\quad\Xi(0,x,\xi)=\xi.
\end{array}
\end{cases}
\end{equation}
It is clear as well that the integral $\int_0^{+\infty} \ldots$ in (\ref{basicform}) carries most of its energy,
semi-classically, over the zero energy level of $h$, defined as
\begin{align}
\label{zeroen}
H_{0}:=\left\{(x,\xi)\in\R^{2d},\ \text{s.t.}\ h(x,\xi)=0\right\}.
\end{align}
In view of the integral  (\ref{basicform}) and of the above considerations, the following definitions are natural.
The first definition is standard.
\begin{def1}
\label{nontra}
{\bf [non-trapping condition]}

\noindent
The refraction index $n^2$ is said {\em non-trapping on the zero energy level}
whenever for each $(x,\xi)\in H_0$, the associated trajectory $(X(t,x,\xi),\Xi(t,x,\xi))$ satisfies
$$
\lim_{t\to+\infty}| X(t,x,\xi)|=+ \infty.
$$
\end{def1}
When the refraction index is non-trapping, the rough idea is that any trajectory $X(t,x,\xi)$ on the zero energy level
leaves any given neighbourhood of the origin $x=0$ in finite time, making the above integral $\int_0^{+\infty} \ldots$ in (\ref{basicform}) converge with respect to the bound $t=+\infty$.

The second definition comes from \cite{MR2139886} (this assumption is called "transversality condition" in the original text).
\begin{def1}
\label{nonre}
{\bf [non-refocusing condition]}

\noindent
We say that $n^2$ satisfies the {\em non-refocusing condition} if the refocusing set,
defined as 
\begin{equation}\label{condrefocusing}
M:=\left\{(t,\xi,\eta)\in]0,+\infty[\times\R^{2d}\ \ s.t.\ \frac{|\eta|^2}{2}=n^2(0),\ X(t,0,\xi)=0,\ \Xi(t,0,\xi)=\eta\right\}
\end{equation}
is such that  $M$ is a submanifold of $]0,+\infty[\times\R^{2d}$ and $M$  satisfies 
$$\text{dim}\, M<d-1.$$
\end{def1}
When the non-refocusing condition is satisfied, the rough idea is that the trajectories $X(t,0,\xi)$ on the zero energy level
issued from the origin $x=0$ at time $t=0$ 
cannot accumulate in any given neighbourhood of the origin $x=0$ at later times $t>0$ (this is encoded in the requirement
on ${\rm dim} \, M$).
Technically speaking, an appropriate stationary phase argument in formula (\ref{basicform}) allows to exploit
in \cite{MR2139886}
the non-refocusing condition and to prove the weak convergence of $w_\eps$ towards~$w_{out}$ under this assumption.
The main result in \cite{MR2139886} is the following: when the refraction index is both non-trapping {\em and}
satisfies the above non-refocusing condition, then  $w_\eps \sim w_{out}$ as $\eps \to 0$ weakly.

Recently, J.F. Bony in \cite{MR2582438} shows the convergence of the {\em Wigner measure} associated with $w_\eps$. He requires a geometrical assumption on the index of refraction that is in the similar spirit, yet weaker, than the above non-refocusing condition,
namely
\begin{equation}
\label{condviriel}
\text{meas}_{n-1}\left\{\xi\in\sqrt{2n^2(0)} \, \Sp^{d-1};\quad \exists\ t>0\quad X(t,0,\xi)=0\right\}=0,
\end{equation}
where $\text{meas}_{n-1}$ is the Euclidian surface measure on $\sqrt{2n^2(0)} \, \Sp^{d-1}$ and $\Sp^{d-1}$ denotes the unit sphere in dimension $d$. Besides, inspired by \cite{MR2139886}, he constructs a  refraction index which is both non-trapping
and  does not satisfy condition (\ref{condviriel}), and in that case he proves the {\em non-uniqueness}
 of the limiting of the Wigner measure.

\medskip

The goal of this paper is to construct a refraction index that is both non-trapping and violates the non-refocusing condition,
and to establish in that case that $w_\eps$ goes weakly to a function of the form
"$w_{out}+$perturbation",  for some explicitly computed and non-zero perturbation.
To be more accurate, we  construct below a refraction index for which the above refocusing manifold 
$\displaystyle  M=\big\{ (t,\xi,\eta)  
\ \text{s.t.} \ \frac{|\eta|^2}{2}=n^2(0),\ X(t,0,\xi)=0,\ \Xi(t,0,\xi)=\eta \big\} $ is smooth, yet has dimension
$\text{dim}\, M=d-1$, a critical case, and we prove
$w_\eps \sim$ "$w_{out}+$perturbation"
in that situation.


\subsection{Construction of the refraction index and statement of our main result}

\begin{figure}[htb]
\begin{center}
\includegraphics[scale=0.6]{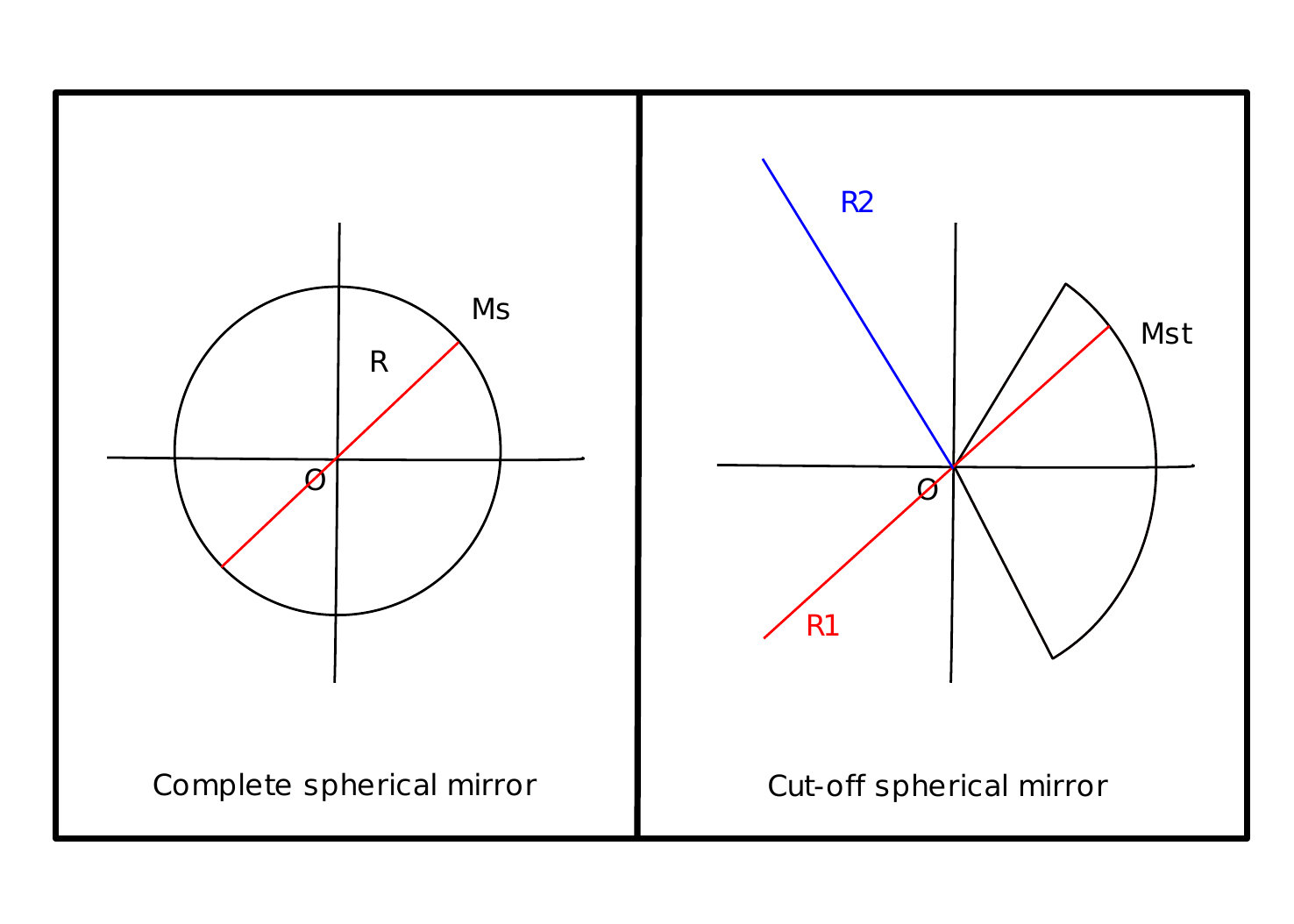}
\caption{Spherical mirror in dimension 2}
\label{fig1}
\end{center}
\end{figure}

Let us first examine the case of  dimension $d=2$. Let $M_{s}$ be a circular mirror centered
at the origin. Any standard ray issued from the origin $x=0$ hits the mirror and goes back to the origin at some later time:
refocusing occurs in a strong fashion.
However all rays are trapped inside the circular mirror, leading to
a trapping situation, in the sense of definition \ref{nontra}. To recover a non-trapping and refocusing situation,
it is necessary
to consider an angular aperture of the circular mirror, with total aperture $<\pi$.
This is shown in figure $\ref{fig1}$: the circular mirror with total aperture $<\pi$ provides a (non-smooth) non-trapping and refocusing refraction index.
To transform the above paradigm into a smooth one,
some regularizations need to be performed. The construction needs to be done in any dimension $d\geq 2$ as well.

Let us first introduce the hyperspherical coordinates $(r,\theta_{1},\ldots,\theta_{d-1})$ in dimension~$d\geq2$
\begin{eqnarray*}
x_{1}&=&r\cos(\theta_{1}),\\
x_{2}&=&r\sin(\theta_{1})\cos(\theta_{2}),\\
x_{3}&=&r\sin(\theta_{1})\sin(\theta_{2})\cos(\theta_{3}),\\
\vdots&&\\
x_{d-1}&=&r\sin(\theta_{1})\ldots\sin(\theta_{d-2})\cos(\theta_{d-1}),\\
x_{d}&=&r\sin(\theta_{1})\ldots\sin(\theta_{d-2})\sin(\theta_{d-1}),
\end{eqnarray*}
with 
$$
\theta_{1}\in [0,\pi], \quad \theta_{j}\in[0,2\pi]  \, \text{ whenever } \, j\geq 2 \text{ when } d\geq 3, \qquad \text{ and }
\theta_1\in [-\pi,\pi] \, \text{ when } d=2.
$$
Next, we choose a {\em fixed}, smooth cut-off function $\chi$ on $\R$ such that
\begin{equation}
\label{chi}
\chi(t)= 1,\quad \forall\ |t|\leq 1,\qquad \chi(t)= 0,\quad \forall\ |t|\geq2,\qquad \chi(t)\geq0,\quad \forall\, t\in\R.
\end{equation} 
We choose a radius $R>0$ and  
define the radial function
\begin{equation}
\label{ff}
f(x)\equiv f(r):=\chi\left(2(r-R)\right),\quad \forall\ x=(r,\theta_{1},\ldots,\theta_{d-1}).
\end{equation}
We choose an angle (aperture) $\theta_{0}\in[0,\pi/4[$,
and  define the angular function
\begin{equation}
\label{gth}g(x)\equiv g(\theta_1):=\chi\left(\frac{\theta_{1}}{\theta_{0}}\right),\quad \forall\ x=(r,\theta_{1},\ldots,\theta_{d-1}).
\end{equation}
a smooth version of the angular aperture $|\theta_1|\leq \theta_0$.
Finally, we choose two parameters $n^2_{\infty}>0$ and $\lambda>0$ such that
\begin{align}
\label{ninf}
 n^2_{\infty}
<
\lambda.
\end{align}
We introduce the following 
\begin{def1}
\label{n22pot}
{\bf [refraction index]}

\noindent
We define the refraction index, retained in the whole subsequent analysis,
as the following smooth version of the circular mirror with total aperture $\theta_0<\pi/4$, namely\footnote{The refraction index is negative in a bounded region of $x$. As already mentioned, we still use the abuse of notation consisting in using the squared of $n$.}
\begin{equation}
\label{n2pot}
n^2(x):=n^2_{\infty}-\lambda f(x)g(x)\equiv n_\infty^2-\lambda f(r) g(\theta_1), \qquad \forall\, x\in\R.
\end{equation}
\end{def1}

\begin{figure}[htb]
\begin{center}
\includegraphics[scale=0.39]{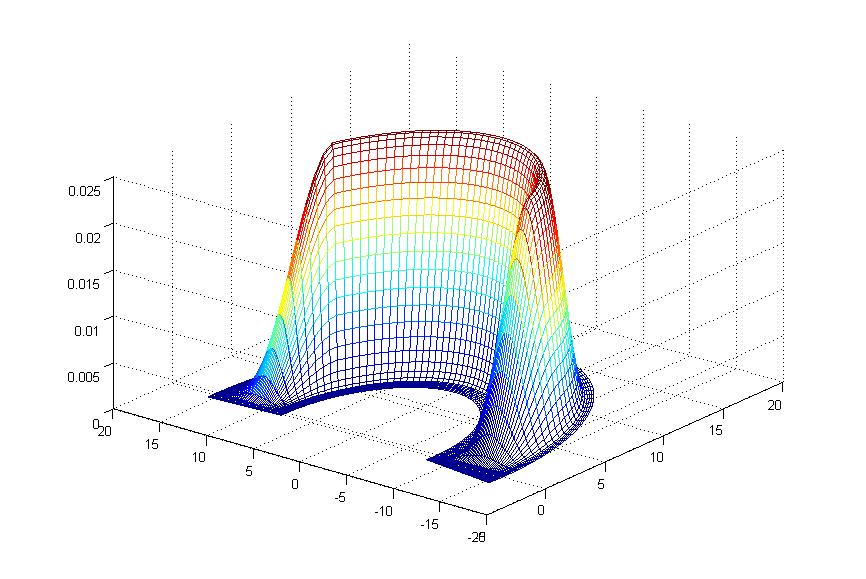}
\caption{The function $n_\infty^2-n^2(x)=\lambda \, f(x) \, g(x)$ in dimension $d=2$}
\label{fig2}
\end{center}
\end{figure}

We are now in position to state our main result. 
Let $(e_{1},\ldots, e_{d})$ be the canonical basis of~$\R^d$. Since the direction $e_1$ is a symmetry axis for our  refraction index,
we introduce for later purposes the space $M_{d}(\R)$  of  square matrices of dimension $d$,
we denote by $\Ot_{d}(\R)$ the space of orthogonal matrices, and we introduce the notation
\begin{equation}
\label{od1}
\Ot_{d,1}(\R):=\left\{ A\in \Ot_{d}(\R),\ \text{s.t.}\ A e_{1}=e_{1}\right\}.
\end{equation}
The refraction index $n^2(x)$ in (\ref{n2pot}) is invariant under the action of  $\Ot_{d,1}(\R)$.
We last introduce a particular set of speeds, namely the set of initial speeds $\xi$ such that the zero energy trajectory
$X(t,0,\xi)$
issued from the origin at time $t=0$  is reflected 
towards the origin at some later time $t>0$. With the retained value of $n^2(x)$, 
we arrive at the definition
\begin{def1}
\label{ith}
{\bf [reflected rays]}

\noindent
The reflection set $I_{\theta_0}$ is defined as
\begin{align}
\label{ith0}
&
\nonumber
I_{\theta_{0}}=\left\{\xi:=(|\xi|,\theta_{1},\ldots,\theta_{d-1}) \in \R^d \, \text{ s.t. } \,
\theta_1\in [-\theta_0,+\theta_0] \, \text{ and } \, \left|\xi\right|=\sqrt{2 \, n^2(0)}\right\}.
\end{align}
\end{def1}
Note that the (intuitive) fact that a velocity $\xi$ is such that $X(t,0,\xi)$ hits the origin at some time $t>0$ if and only if
$\xi \in I_{\theta_0}$, is proved later (see section \ref{statset}).

\medskip

Our main result in this text is the
\begin{theo}
\label{propnonconvergence}
{\bf [Main Result]}

\noindent
Let $n^2$ be the refraction index defined in  (\ref{n2pot}).
Assume the aperture $\theta_0<\pi/4$ and the radius $R>0$ satisfy the smallness condition  
\begin{equation}
\label{conditionpetitpotentiel}
1-\cos(2\theta_0) < \frac{1}{2 R}.
\end{equation}
Assume $d\geq 3$. Then, the following holds:\\

\noindent
\textit{i)} The index $n^2$ is non-trapping on the zero-energy level $H_{0}=\{(x,\xi) \, \text{ s.t. } \, |\xi|^2/2-n^2(x)=0\}$.
\\
\\
\textit{ii)} The refocusing set $M=\{(t,\xi,\eta) \, \text{ s.t. } \, |\eta|^2=2 n^2(0), \,
X(t,0,\xi)=0, \,
\Xi(t,0,\xi)=\eta\}$ (see (\ref{condrefocusing}))
is a smooth submanifold of $]0,+\infty[\times \R^{2d}$, with boundary, and its dimension has the critical value
$$
{\rm dim}(M)=d-1.
$$
\textit{iii)} Assume
 the source term $S$ satisfies $S\in\Sc(\R^d)$.
%
%
Then, we have
\begin{equation*}
\forall\,\phi\in\Sc(\R^d),\qquad
\big\langle
w_{\varepsilon}-\left( w_{out}+L_\eps \right) \, , \, 
\phi  \big\rangle  \mathop{\longrightarrow}_{\eps\to 0} 0,
\end{equation*} 
where the distribution $L_\eps$ is defined for any $\phi\in\Sc(\R^d)$ through
\begin{align}
\label{leps}
&
\nonumber
\langle L_{\varepsilon }, \phi \rangle= \\
&
\quad
C_{n^2,d} 
\int_{I_{\theta_{0}}} 
\exp\left(
\frac{i}{\varepsilon} \int_{0}^{T_{R}} \left(
\frac{\left|\Xi(s,0,\xi)\right|^2}{2}
+
n^2(X(s,0,\xi))
\right)  \, ds \right) \, 
\widehat S(\xi)\widehat{\phi^*}(-\xi) \,  d\sigma_{\theta_0}(\xi).
\end{align} 
Here $d\sigma_{\theta_0}$ denote the natural Euclidean surface measure on $I_{\theta_0}$ (see definition \ref{ith}),
the return time $T_R>0$ is the unique time\footnote{The fact that all these quantities exist and are well defined
is part of the Theorem, and is proved in section \ref{statset}.}
such that for any $\xi \in I_{\theta_0}$ we have $X(T_R,0,\xi)=0$,
and the constant $C_{n^2,d}\neq 0$ can be explicitly computed and depends  only on the index $n^2$ and on the
dimension~$d$.
\end{theo}

\begin{rem}
The condition \eqref{conditionpetitpotentiel} is technical, and requires the aperture $\theta_0$ to be small:
it
ensures  the trajectories cannot be trapped by the refraction index.
\end{rem}

\begin{rem}
Note in passing that the constraint $d\geq 3$, which is also needed in  reference \cite{MR2139886},
comes from a stationary phase argument. This constraint on the dimension is standard in the analysis of
Schr\"odinger-like operators. It comes from the fact that the dispersion induced by the free Schr\"odinger operator acts like $t^{-d/2}$,
a factor that is integrable close to $t=+\infty$ whenever $d\geq 3$.
\end{rem}

\begin{rem}
Let $\xi_0:=(\sqrt{2} \, n(0) , 0 ,\ldots,0)$.
The distribution $L_\eps$ can as well be written as
\begin{align*}
&
\nonumber
\langle L_{\varepsilon} \, ,  \, \phi \rangle= \\
&
\quad
C_{n^2,d} 
\,
\exp\left(
\frac{i}{\varepsilon} \int_{0}^{T_{R}} \left(
\frac{\left|\Xi(s,0,\xi_0)\right|^2}{2}
+
n^2(X(s,0,\xi_0))
\right)  \, ds \right) \, 
\left(
\int_{I_{\theta_{0}}} 
\widehat S(\xi)\widehat{\phi^*}(-\xi) \, d\sigma_{\theta_0}(\xi)
\right).
\end{align*} 
This formulation illustrates in a clearer way the fact
that if the source $S$ radiates towards the mirror, then $w_{\varepsilon}$ converges
towards a non-trivial perturbation of $w_{out}$.

Note in passing that in the present counter-example, as in the paper by J.-F. Bony
\cite{MR2582438} , only {\em subsequences} of $w_\eps$ converge, due to the above oscillatory factor
$\exp(i \,  {\rm const.}/\eps)$. 
\end{rem}

\begin{rem}
In the chosen hyperspherical coordinates, the Euclidean measure
$d\sigma_{\theta_0}(\xi)$ coincides with
$\displaystyle
d\sigma_{\theta_0}(\xi)= n(0)^{d-1} \, d\sigma(\theta_1, \ldots, \theta_{d-1})$, 
where
$d\sigma(\theta_1, \ldots, \theta_{d-1})$ denotes the standard euclidean surface measure on
the unit sphere ${\mathbb S}^{d-1}$.
\end{rem}


\subsection{Preliminary reduction of the proof}


Our main result contains three distinct statements.
Items (i) and (ii) are of geometric nature, and merely concern the behavious of the classical trajectories associated with 
the retained refraction index. Their proof is performed in sections \ref{behaviour} and \ref{statset}, respectively.
Item (iii) is the main item, and concerns the asymptotic analysis of $w_\eps$.
Since our analysis heavily relies on tools previously developped in  \cite{MR2139886}, we briefly recall here some of these
tools and indicate how the analysis of $w_\eps$ can be reduced to a simpler sub-problem. We postpone the analysis of the
reduced subproblem, hence of item (iii) of our main result,  to section \ref{proposition} below.

\medskip

As already indicated, given a smooth test function $\phi$, we start from the formulation
\begin{equation*}
\left\langle w_{\varepsilon},\phi\right\rangle=\frac{i}{\varepsilon}\int_{0}^{+\infty} e^{-\alpha_{\varepsilon}t}
\big\langle U_{\varepsilon}(t)S_{\varepsilon},\phi_{\varepsilon}\big\rangle \,
dt.
\end{equation*}
(See above for the notation). 
The next step consists in splitting the above time integral into four time scales, namely
very small, small, moderate, and large time scales. To do so, we take one small parameter $\theta>0$ and two large parameters
$T_0>0$ and $T_1>0$, and split the above time integral into the four zones
$$
0 \leq t \leq T_0 \, \eps, \quad
T_0 \, \eps \leq t \leq \theta, \quad
\theta \leq t \leq T_1, \quad
T_1\leq t \leq +\infty \quad
(\theta \ll 1, \quad T_0, T_1 \gg 1).
$$
Technically, we use a smooth splitting, based on the already used cut-off function $\chi$ (see (\ref{chi})).
Besides, we also distinguish between the contribution of zero and non-zero energies, namely taking a small parameter $\delta>0$, we write, in the sense of functional caculus for self-adjoint operators, the identity
$$
1=\chi_\delta(H_\eps) + (1-\chi_\delta)(H_\eps), \, \text{ where } \, H_\eps:=\frac{\eps^2}{2} \Delta_x  + n^2(x),
\, \text{ and } \, \chi_\delta(s):=\chi\left(\frac{s}{\delta} \right) \, (s \in \R, \delta\ll 1).
$$

The main intermediate result of the present subsection is the following
\begin{prop}
\label{intermres}
{\bf [Main intermediate result]}

\noindent
Take a test function $\phi\in\Sc(\R^d)$.
Define $\widetilde{w_\eps}$ as
\begin{equation*}
\left\langle \widetilde{w_{\varepsilon}},\phi\right\rangle:=
\frac{i}{\varepsilon}
\int_{\theta}^{T_1} 
\left( 1 - \chi\right)\left( \frac{t}{\theta}\right)  \,
e^{-\alpha_{\varepsilon}t} \, 
\big\langle U_{\varepsilon}(t) \, 
\chi_\delta(H_\eps)  \,
S_{\varepsilon},\phi_{\varepsilon}\big\rangle \, dt.
\end{equation*}
Then,
there is a large $T_1>0$ such that for any small  $\delta>0$,  and any small   $\theta>0$, 
there exists a constant
$C_{\theta,\delta}>0$ such that for any small $\eps>0$, we have
\begin{equation*}
\Big|
\Big\langle \left( w_\eps - (w_{out} + \widetilde{w_{\varepsilon}}) \right) \, , \, \phi\Big\rangle
\Big|
\leq
C_{\theta,\delta} \, \left(
\frac{1}{T_0^{d/2-1}}
+
\frac{1}{T_0} + \a_\eps^2
+
\eps
\right).
\end{equation*}
\end{prop}

This result roughly asserts that $w_\eps$ is asymptotic to  $w_{out}+\widetilde{w_\eps}$ as $\eps\to 0$,
up to carefully choosing the various parameters $T_0$, $T_1$, etc. Hence the proof of item (iii) of our main result
essentially reduces to proving that $\widetilde{w_\eps}\sim L_\eps$ as $\eps\to 0$.

\medskip

\noindent
{\bf Proof of Proposition \ref{intermres}.}

The proof is obtained by gathering the statements of Proposition
\ref{propvst}, Proposition
\ref{propoutH}, Proposition 
\ref{propHlt}, Proposition 
\ref{propdur} below.
\vspace{-0,4cm}
\hfill$\blacksquare$\\

\medskip

The remainder part of this paragraph is devoted to a brief idea of proof of the above 
auxiliary Propositions that lead to Proposition \ref{intermres}.

\medskip

\paragraph{$\bullet$ Contribution of very small times $0 \leq t \leq T_0 \, \eps$.}

\mbox{}\\

The contribution of very small times to 
$\left\langle w_{\varepsilon},\phi\right\rangle=i \, \eps^{-1} \displaystyle \, \int_{0}^{+\infty} e^{-\alpha_{\varepsilon}t}
\left\langle U_{\varepsilon}(t)S_{\varepsilon},\phi_{\varepsilon}\right\rangle \,
dt$, is
$$\frac{i}{\varepsilon}\int_{0}^{2T_{0}\varepsilon}\chi\left(\frac{t}{T_{0}\varepsilon}\right)e^{-\alpha_{\varepsilon}t}
\big\langle U_{\varepsilon}(t)S_{\varepsilon},\phi_{\varepsilon}\big\rangle \, dt.$$
It is the main contribution to $w_\eps$, provided $T_{0}$ is large enough. Indeed, we have the following fact, whose proof is based on a simple weak convergence argument.
\begin{prop} 
\label{propvst}
{\bf (See  \cite{MR2139886}).} \,
Let $n^2(x)$ be {\em any} bounded and continuous refraction index. Then, if S and $\phi$ belong to  $\Sc(\R^d)$, we have

\medskip

(i) For all time $T_0>0$,
\begin{equation*}
\frac{i}{\varepsilon}\int_0^{2T_0\varepsilon} \chi\left(\frac{t}{T_0\varepsilon}\right)e^{-\alpha_\varepsilon t}
\big\langle U_\varepsilon(t)S_\varepsilon,\phi_\varepsilon\big\rangle \, dt\underset{ \varepsilon\to 0}{\longrightarrow}i\int_0^{2T_0} \chi\left(\frac{t}{T_0}\right)\left\langle \exp\left(it\left(\frac{\Delta_x}{2}+n^2(0)\right)\right)S,\phi\right\rangle \, dt.
\end{equation*}

\medskip

(ii)
There exists $C_d>0$ which only depends on the dimension such that
\begin{equation*}
\left|
\left(
\frac{i}{\varepsilon}\int_0^{2T_0\varepsilon} \chi\left(\frac{t}{T_0}\right)
\big\langle \exp(it(\Delta_x/2+n^2(0)))S,\phi\big\rangle \, dt
\right)-
\left\langle w_{out},\phi\right\rangle\right|\\
\leq \frac{C_d}{T_0^{d/2-1}}.
\end{equation*}
\end{prop}

\paragraph{$\bullet$ Contribution of small, up to large times, away from the zero-energy level.}

\mbox{}\\

The contribution  to
$\left\langle w_{\varepsilon},\phi\right\rangle=i \, \eps^{-1} \displaystyle \, \int_{0}^{+\infty} e^{-\alpha_{\varepsilon}t}
\big\langle U_{\varepsilon}(t)S_{\varepsilon},\phi_{\varepsilon}\big\rangle \,
dt$ that is associated with small, up to large times, away from the zero-energy level, is
$$\frac{i}{\varepsilon}\int_{T_{0}\varepsilon}^{+\infty} e^{-\alpha_{\varepsilon} t}(1-\chi)\left(\frac{t}{T_{0}{\varepsilon}}\right)
\big\langle (1-\chi_{\delta})(H_{\varepsilon})U_{\varepsilon}(t)S_{\varepsilon},\phi_{\varepsilon}\big\rangle \, dt.$$
It is seen to be small, using a non-stationary phase argument in time, see \cite{MR2139886} (this is the reason for the previous
cut-off close to the initial time $t=0$, where integrations by parts in time are forbidden). Indeed, we have the
\begin{prop}
\label{propoutH}
{\bf (See  \cite{MR2139886}).} \,
Let $n^2$ be {\em any} long-range refraction index. Let $S$ and $\phi$ belong to $L^2(\R^d)$. Then there exists a constant $C_\delta>0$, which only depends on $\delta>0$, such that for any small $\varepsilon>0$ and any  $T_0>0$,
we have
\begin{eqnarray*}
\left|\frac{1}{\varepsilon}\int_{T_0\varepsilon}^{+\infty}(1-\chi)\left(\frac{t}{T_0\varepsilon}\right)\big\langle\left(1-\chi_\delta(H_\varepsilon)\right)U_\varepsilon(t)S_\varepsilon,\phi_\varepsilon(t)\big\rangle \, dt\right|\leq C_\delta\left(\frac{1}{T_0}+\alpha_\varepsilon^2\right).
\end{eqnarray*}
\end{prop}

\paragraph{$\bullet$ Contribution of large times, near the zero-energy level.}

\mbox{}\\

The contribution  to
$\left\langle w_{\varepsilon},\phi\right\rangle=i \, \eps^{-1} \displaystyle \, \int_{0}^{+\infty} e^{-\alpha_{\varepsilon}t}
\big\langle U_{\varepsilon}(t)S_{\varepsilon},\phi_{\varepsilon}\big\rangle \,
dt$ that is associated with  large times, close to the zero-energy level, is
$$\frac{i}{\varepsilon}\int_{T_{1}}^{+\infty} e^{-\alpha_{\varepsilon}t}\left\langle \chi_{\delta}(H_{\varepsilon})U_{\varepsilon}(t)S_{\varepsilon},\phi_{\varepsilon}\right\rangle dt.$$
It is seen to be of order $O(\varepsilon^N)$, for all $N\in\N$, see \cite{MR2139886}. Indeed,
the semiclassical support of $\chi_{\delta}(H_{\varepsilon})U_{\varepsilon}(t)S_{\varepsilon}$ 
goes to infinity in the $x$ direction at speed of the order $1$ ({\em i.e.} the semi-classical support
lies in a region that is at distance of order $t$ from the origin -- this uses an argument due to Wang, see
\cite{MR927007}), while the semi-classical support of $\phi_{\varepsilon}$ remains close to the origin.
 This argument relies on the fact that for $T_{1}$ large enough, the semiclassical supports of the two functions are disconnected, which in turn uses the non-trapping behaviour of the refraction index.
We arrive at
\begin{prop}
\label{propHlt}
{\bf (See  \cite{MR2139886}).} \,
 Let $n^2$ be {\em any} long-range refraction index that is non-trapping. 
Let S and $\phi$ be in $\Sc(\R^d)$. Then there exist $\delta_0>0$ and $T_1(\delta_0)>0$ such that for all time $T_1\geq T_1(\delta_0)$ and any $0<\delta<\delta_0$, there exists a constant $C_{\delta}$ such that 
\begin{eqnarray*}
\left|\frac{1}{\varepsilon}\int_{T_1}^{+\infty}
e^{-\alpha_\varepsilon t}
\big\langle\chi_\delta(H_\varepsilon)U_\varepsilon(t)S_\varepsilon,\phi_\varepsilon\big\rangle \, dt \right|
&\leq& C_{\delta}\varepsilon.
\end{eqnarray*}
\end{prop}

\paragraph{$\bullet$ Contribution of small times near the zero-energy level}

\mbox{}\\

The contribution  to
$\left\langle w_{\varepsilon},\phi\right\rangle=i \, \eps^{-1} \displaystyle \, \int_{0}^{+\infty} e^{-\alpha_{\varepsilon}t}
\big\langle U_{\varepsilon}(t)S_{\varepsilon},\phi_{\varepsilon}\big\rangle \,
dt$ that is associated with  small times, close to the zero-energy level, is
$$
\frac{i}{\varepsilon} \int_{T_0 \eps}^{\theta}
e^{-\alpha_{\varepsilon}t} \, (1-\chi)\left(\frac{t}{T_0 \eps}\right)
\chi\left(\frac{t}{\theta}\right) \,
\big\langle U_{\varepsilon}(t)\chi_{\delta}(H_{\varepsilon})S_{\varepsilon},\phi_{\varepsilon}\big\rangle dt.
$$
Unlike in the previous case, the semiclassical supports of $U_{\varepsilon}(t)\chi_{\delta}(H_{\varepsilon})S_{\varepsilon}$
 and $\phi_{\varepsilon}$ {\em may} intersect for these values of time $t$. 
The whole point in \cite{MR2139886} lies, roughly speaking,
in proving a {\em dispersion estimate}. The key is to prove that the variable coefficients
Schr\"odinger propagator $U_\eps(t)$ has the same dispersive properties than the free Schr\"odinger propagator,
corresponding to the case when $n^2\equiv 0$,
at least for small values of $t$ such that $0 \leq  t \leq \theta$ (for later times, the semiclassical support
of $U_\eps(t) \, S_\eps$ is close to the classical trajectories $(X(t),\Xi(t))$, trajectories
which in turn may come back close to the origin and contradict any dispersion effect).
Indeed, for small times, the trajectory $(X(t),\Xi(t))$ is close to its first order expansion in time,
which is the key to obtaining dispersive effects similar to the one
at hand in the free case. Technically speaking, the proof relies on establishing that the propagator
$U_\eps(t)$ behaves like the free Schr\"odinger propagator for small times,
a propagator whose symbol is $\exp(i t |\xi|^2/\eps)$, and
which in turn has size $(\eps/t)^{d/2}$ thanks to a stationary phase argument.

To obtain the desired statement, a wave packet approach is actually introduced, which strongly uses
the work by Combescure and Robert $(\cite{MR1461126})$. It allows to compute explicitly the propagator $U_{\varepsilon}(t) \, S_{\varepsilon}$, using the Hamiltonian flow and related, linearized, quantities,
to obtain a representation of the form
\begin{align}
\label{eqJeps}
\nonumber
&
\frac{i}{\varepsilon} \int_{T_0 \eps}^{\theta}
e^{-\alpha_{\varepsilon}t} \, (1-\chi)\left(\frac{t}{T_0 \eps}\right)
\chi\left(\frac{t}{\theta}\right) \,
\big\langle U_{\varepsilon}(t)\chi_{\delta}(H_{\varepsilon})S_{\varepsilon},\phi_{\varepsilon}\big\rangle  \, dt
\\
&
\qquad\qquad
=
\frac{1}{\varepsilon^{(5d+2)/2}} \,
\int_{T_{0}\varepsilon}^\theta \, \int_{\R^{6d}}
e^{\frac{i}{\varepsilon} \, \psi(t,X)}  \, a_{N}(t,X)
\, dt \, dX
+ O_{\theta,\delta}\left(\eps^N\right),
\end{align}
where $X=(q,p,x,y,\xi,\eta) \in \R^{6d}$, where $N$ is a possibly large integer,
and the remainder term $O_{\theta,\delta}\left(\eps^N\right)$ is upper bounded by $C_{\theta,\delta} \, \eps^N$ for some $C_{\theta,\delta}>0$ independent of $\eps$, which depends on the chosen $\theta>0$ and $\delta>0$.
Note  that the amplitude $a_{N}$  is defined in (\ref{amplitude}) below,
while the {\em complex} phase function $\psi$ is defined in (\ref{phase}) below.
We refer to section \ref{proposition} for details about the representation formula
(\ref{eqJeps}), which is a key ingredient in our proof of the main theorem.

\medskip

With this representation at hand, we arrive at the
\begin{prop}
\label{propdur}
{\bf (See  \cite{MR2139886}).} \,
Let $n^2$ be {\em any} long-range potential which is non-trapping. For  $\theta$ and $\delta$ small enough,
there exists $C_\theta>0$ and $C_{\theta,\delta}>0$ such that for all $\varepsilon\leq 1$ we have
\begin{align}
\frac{1}{\varepsilon}\int_{T_0\varepsilon}^{\theta}\chi\left(\frac{t}{\theta}\right)\left(1-\chi\left(\frac{t}{T_0\varepsilon}\right)\right)e^{-\alpha_\varepsilon t}\left\langle U_\varepsilon(t)\chi_\delta(H_\varepsilon)S_\varepsilon,\phi_\varepsilon\right\rangle dt
\leq \frac{C_\theta}{T_0^{d/2-1}}+C_{\theta,\delta} \, \eps.
\label{propfinal2}
\end{align}
\end{prop}


\section{Properties of the refraction index}
\label{Sectionconstruction}



\subsection{Non-trapping behaviour}\label{behaviour}


The goal of this subsection is to prove item (i) of our main Theorem \ref{propnonconvergence}.

We prove that the chosen refraction index $n^2(x)=n^2_{\infty}-\lambda f(r)g(\theta_{1})$
in (\ref{n2pot}) is \textit{non-trapping} on the zero-energy level
$H_{0}=\left\{(x,\xi)\in\R^{2d},\ \text{s.t.}\ \xi^2/2=n^2(x)\right\}$.

\medskip

We first observe that the zero energy level has the more explicit value
\begin{align*}
H_{0}
=\left\{(x,\xi)\in\R^{2d},\ \text{s.t.}\ x=(r,\theta_{1},\ldots,\theta_{d-1}),\ \frac{\xi^2}{2}=n^2_{\infty}-\lambda f(r)g(\theta_{1})\right\}.
\end{align*}
We readily define the following two regions. The first one is usually called the classically forbidden region:
any trajectory living on the zero-energy level cannot reach the set $B_\emptyset$. The second one is sometimes called  here the bump of the refraction index: it is the region where the refraction index actually {\em varies} with $x$.
Outside this region, the refraction index is constant and the Hamiltonian trajectories associated with
$h(x,\xi)=|\xi|^2/2+n^2(x)$ are straight lines.
\begin{def1}
\label{defini}
(i)
We denote by $B_{\emptyset}$ the set (classically forbidden region)
\begin{equation*}
B_{\emptyset}:=\left\{x\in\R^d, \ \text{s.t.}\ n^2(x)<0\right\}=\left\{x=(r,\theta_{1},\ldots,\theta_{d-1}),\ \text{s.t.}\ n^2_{\infty}<\lambda f(r)g(\theta_{1}) \right\}.
\end{equation*}
(ii)
We denote by $B_p$ the set  (bump)
\begin{align*}
B_{p}:=
\left\{x=(r,\theta_{1},,\ldots,\theta_{d-1}),\ \text{s.t.}\ R-1\leq r\leq R+1,\,|\theta_{1}| \leq2\theta_{0}\right\}.
\end{align*}
\end{def1}

\begin{rem}
From the definition of $B_\emptyset$ and the two functions $f(r)=\chi(2(r-R))$ and $g(\theta_1)=\chi(\theta_1/\theta_0)$
it is clear that there exists $\mu \in ]1,2[$ such that
\begin{align}
\label{bvide}
B_\emptyset
\subset
\left\{R-\frac{\mu}{2}\leq r \leq R+\frac{\mu}{2}, \, |\theta_1|\leq \mu \theta_0\right\}.
\end{align}
It suffices to take $\mu$ such that
\begin{align}
\label{mu}
0
<
\chi(\mu)
<
\frac{n_\infty}{\sqrt{\lambda}}
\, \text{ (hence  } \mu \in ]1,2[\text{)}.
\end{align}
\end{rem}

 \begin{figure}[htb]
\begin{center}
\includegraphics[scale=0.5]{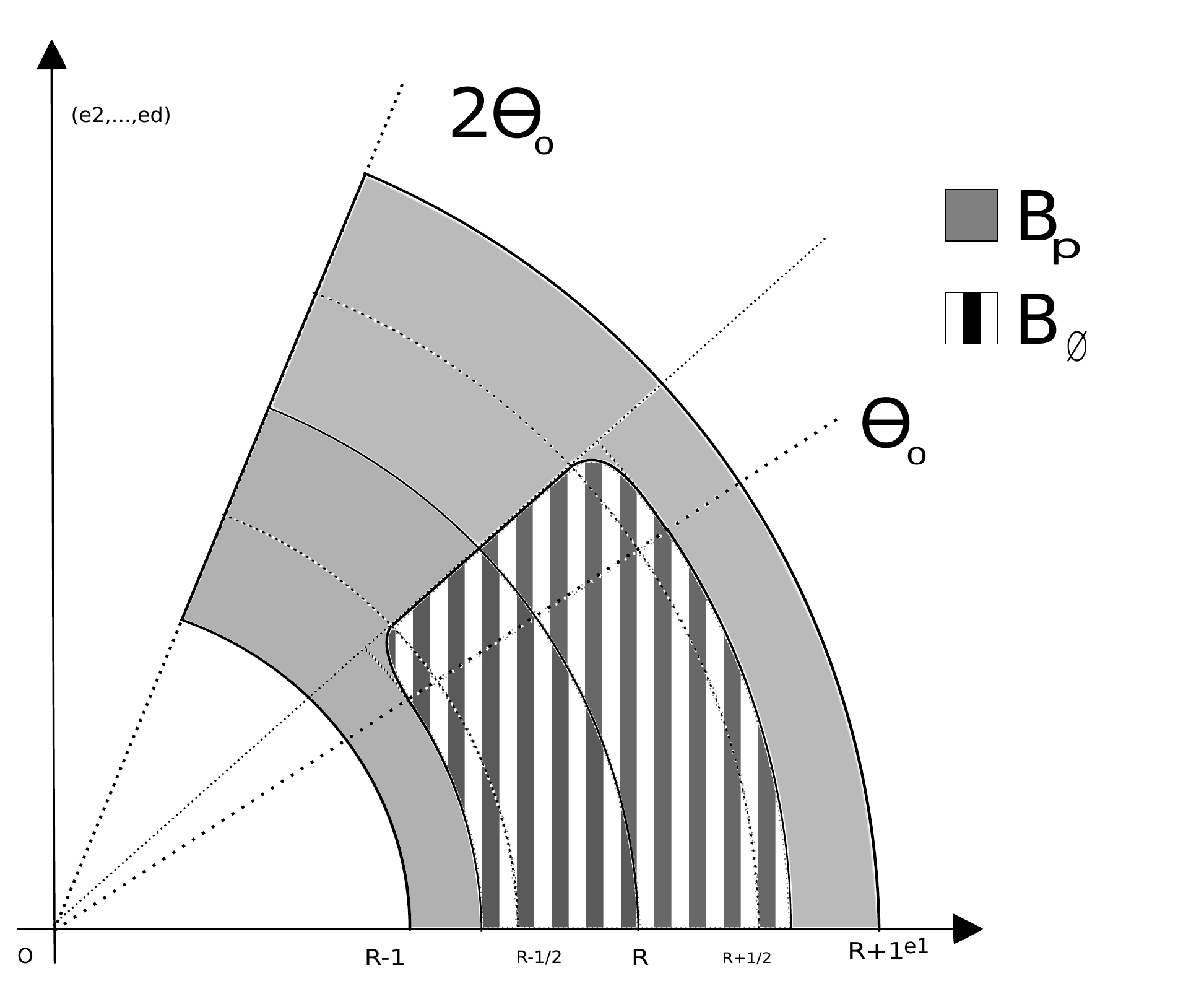}
\caption{Bump of refraction index, and classically forbiden region}
\label{fig3}
\end{center}
\end{figure}

\medskip

Our main step lies in proving the following 
escape estimate
\begin{lem}\label{propasympt}
Select the refraction index $n^2(x)$ as in (\ref{n2pot}) and
assume condition \eqref{conditionpetitpotentiel} is fulfilled, namely $1-\cos(2 \theta_0)<1/(2R)$.
Take a Hamiltonian trajectory $X(t,x,\xi)\equiv X(t)$ living on the zero-energy level and define $x_{0}:=(R,0,\ldots,0)$
in Cartesian coordinates. 

Then, there exists $\alpha>0$, as well as $\beta\in\R$ and $\gamma\in \R$, such that
\begin{equation*}
\forall \,t\geq 0,\quad \left|X(t)-x_{0}\right|^2\geq \alpha\, t^2+\beta\, t+\gamma.
\end{equation*}
\end{lem}
An immediate corollary of the above Lemma  is
\begin{cor}
\label{cornontrapping}
Assume condition \eqref{conditionpetitpotentiel} is fulfilled, namely $1-\cos(2 \theta_0)<1/(2R)$.
Then the refraction index $n^2(x)$ in (\ref{n2pot}) is non-trapping  on the zero-energy level.
\end{cor}

\noindent
{\bf Proof of Corollary \ref{cornontrapping}.}

Apply the preceding lemma and  let $t\to +\infty$.
\vspace{-0,4cm}
\hfill$\blacksquare$\\

\medskip

\noindent
{\bf Proof of Lemma \ref{propasympt}.}
\nopagebreak
\nopagebreak
\paragraph{$\bullet$ First step.} 
We compute the second derivative of $\left|X(t)-x_{0}\right|^2$ and get
\begin{align}
\label{eqsecondderivative}
&
\nonumber\frac{1}{2}\frac{d^2}{dt^2}\left|X(t)-x_{0}\right|^2
=
\left\langle\frac{d^2}{dt^2}X(t),X(t)-x_{0}\right\rangle+\left|\frac{d X}{dt}(t)\right|^2,
\\
\nonumber
&
\qquad
=
\left\langle\nabla n^2(X(t)),X(t)-x_{0}\right\rangle+\left|\frac{d X}{dt}(t)\right|^2,
\\
&
\qquad
=
\left\langle\nabla n^2(X(t)),X(t)-x_{0}\right\rangle+n^2(X(t)),
\end{align}
where we have used the fact that the Hamiltonian trajectory $(X(t),\Xi(t))$ belongs to $H_{0}$.
Letting $X(t)=r\,\vec{u}_{r}$ in hyperspherical coordinates and $x_{0}=(R,0,\ldots,0)$ in Cartesian coordinates,
we obtain on the other hand
\begin{align*}
\left\langle\nabla n^2(X(t)),X(t)-x_{0}\right\rangle&=\left\langle-\lambda f'(r)g(\theta_{1})\vec{u}_{r}-\lambda\frac{f(r)}{r}g'(\theta_{1})\vec{u}_{\theta_{1}},r\,\vec{u}_{r}-R\,\vec{e}_{1}\right\rangle,\\
&=F_{r}(r,\theta_{1})+F_{\theta}(r,\theta_{1}),
\end{align*}
where
\begin{align}
\label{frft}
F_{r}(r,\theta_{1})=-\lambda  f'(r)g(\theta_{1})\left(r-R\cos(\theta_{1})\right),\qquad
F_{\theta}(r,\theta_{1})=-\lambda\frac{R}{r}f(r) g'(\theta_{1}) \sin(\theta_{1}).
\end{align}
Eventually we have
\begin{equation}
\label{npo}
\frac{1}{2}\frac{d^2}{dt^2}\left|X(t)-x_{0}\right|^2=F_{r}(r,\theta_{1})+F_{\theta}(r,\theta_{1})+n^2(X(t)).
\end{equation}
Therefore, the lemma is proved once we establish
the existence of $\alpha>0$ such that
$$
F_{r}(x)+F_{\theta}(x)+n^2(x)\geq \alpha>0
$$
whenever $x \in \Pi_x H_0=\R^d\setminus B_\emptyset$ (where
$\Pi_x$ denotes the projection $(x,\xi)\mapsto x$ from $\R^{2d}$ to $\R^d$).

We readily notice that $n^2$ and $F_{\theta}$ are clearly non-negative function on the whole of $\R^d$.

\paragraph{$\bullet$ Step two: non-negativity of $F_{r}$.}
First, on  $\R^d\setminus B_{p}$, the function $F_{r}$ is zero, hence non-negative. In the same way on  
$B_p\cap \{R-1/2\leq r\leq R+1/2\}$, we have $f'\equiv0$, hence $F_{r}\equiv0\geq0$. There remains to study the
non-negativity of $F_{r}$ on the two sets  $\left\{R-1\leq r\leq R-1/2, \, |\theta_{1}|\leq2\theta_{0}\right\}$
and $\left\{R+1/2\leq r\leq R+1, \, |\theta_{1}|\leq2\theta_{0}\right\}$.

On $\left\{R-1\leq r\leq R-1/2, \, |\theta_{1}|\leq2\theta_{0}\right\}$, we have
\begin{align*}
r-R\cos(\theta_{1})&\leq R-\frac{1}{2}-R\cos(2\theta_{0})= R(1-\cos(2\theta_{0}))-\frac{1}{2}< 0,
\end{align*}
thanks to our assumption (\ref{conditionpetitpotentiel}).
Since $f'\geq0$ on $\left\{R-1\leq r\leq R-1/2\right\}$, we get
$F_{r}\geq 0$ on $\left\{ R-1\leq r\leq R-1/2, \, |\theta_{1}|\leq2\theta_{0}\right\}$.
A similar computation proves that 
$F_{r}\geq 0$ on the set $\left\{ R+1/2\leq r\leq R+1, \, |\theta_{1}|\leq2\theta_{0}\right\}$.

We have obtained that $F_{r} \geq 0$ on the whole of $\R^d$.

\paragraph{$\bullet$ Step three: decomposition of $\R^d$.}
We have just proved that $F_r(x)+F_\theta(x)+n^2(x) \geq 0$ for all $x \in \R^d$.
We now wish to obtain a positive lower bound for $x\notin B_\emptyset$.
The argument relies on the fact that  the refraction index $n^2$ is positive away from the boundary
$\partial B_{\emptyset}$, where
$\partial B_{\emptyset}:=\{(r,\theta_{1},\ldots,\theta_{d-1}),\quad f(r)g(\theta_{1})=n^2_{\infty}/\lambda\}$, while
the term $F_r+F_\theta$ stemming from the gradient of the refraction index in (\ref{npo}) is positive close to
the boundary 
$\partial B_\emptyset$. This is the reason for the decomposition we now introduce.

We define the set (piece of ring)
\begin{align*}
C_{\alpha,\beta}:=\left\{ R-\alpha\leq r\leq R+\alpha,\ -\beta\leq\theta_{1}\leq\beta\right\}.
\end{align*}
We know from the remark after Definition \ref{defini} that there exist $\mu \in ]1,2[$
such that
$$
B_{\emptyset}\subset C_{R+\mu/2,\mu \theta_0}.
$$
We therefore decompose 
$$
\R^d\setminus B_\emptyset
= \left(\R^d\setminus C_{R+\mu/2,\mu \theta_0}\right) \cup  \left(C_{R+\mu/2,\mu \theta_0}\setminus B_\emptyset \right).
$$
We readily observe that,
by construction of $\mu$ (namely $\chi(\mu)^2\in ]0,n_\infty^2/\lambda[$ -- see (\ref{mu})),  for any
$x\in \R^d\setminus C_{R+\mu/2,\mu \theta_0}$, we have the lower bound
$$
n^2(x)=n_\infty^2-\lambda f(r) g(\theta_1)
\geq n_\infty^2-\lambda \chi(\mu)^2
=:c_{n^2}
>0,
$$
There only remains to prove the existence of $c_\nabla>0$ such that $F_r+F_\theta\geq c_\nabla$ on
$C_{R+\mu/2,\mu \theta_0}\setminus B_\emptyset$.

\paragraph{$\bullet$ Step four: positive lower bound for $F_r+F_\theta$ on $C_{R+\mu/2,\mu \theta_0}
\setminus B_\emptyset$.}
Take  $\nu \in ]1,2[$ such that
$$
\frac{n_\infty}{\sqrt{\lambda}} < \chi(\nu) < 1.
$$
where $\chi$ is the truncation function defined in (\ref{chi}).
With this choice of $\nu$, we clearly have,
whenever $x \in C_{R+\nu/2,\nu \theta_0}$, the relation
$n^2(x)=n_\infty^2-\lambda \chi(2(r-R)) \chi(\theta_1/\theta_0)\leq
n_\infty^2-\lambda \chi(\nu)^2<0$, hence
\begin{align*}
C_{R+\nu/2,\nu \theta_0}
\subset B_\emptyset
\subset
C_{R+\mu/2,\mu \theta_0}.
\end{align*}
Therefore, it is enough to obtain a lower bound on $F_r+F_\theta$ on the set
$C_{R+\mu/2,\mu \theta_0}\setminus C_{R+\nu/2,\nu \theta_0}$.

\begin{figure}[htb]
\begin{center}
\includegraphics[scale=0.6]{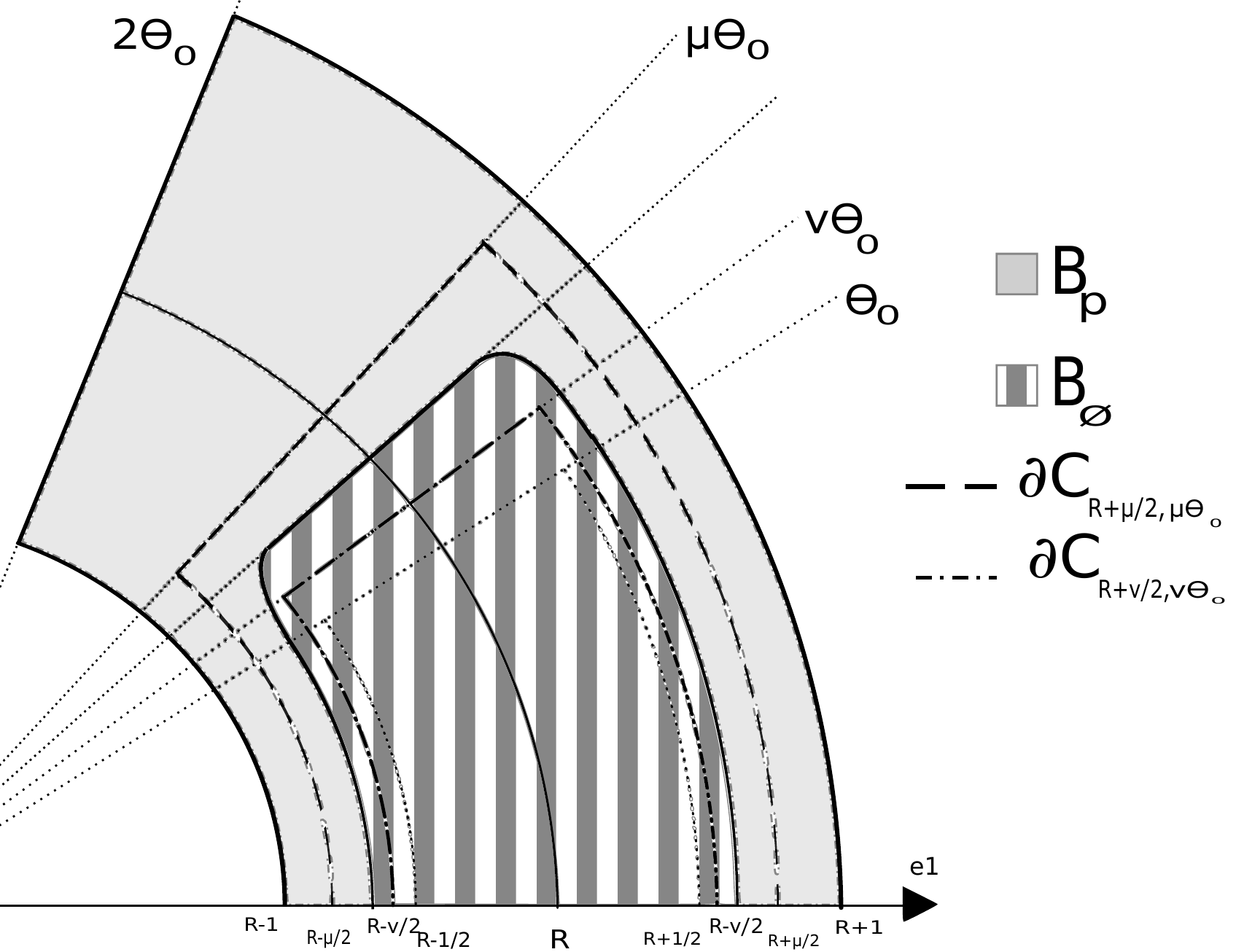}
\caption{Zone of study}
\label{fig4}
\end{center}
\end{figure}

To this end, we decompose (see Figure \ref{fig4})
\begin{align*}
&
C_{R+\mu/2,\mu \theta_0}\setminus C_{R+\nu/2,\nu \theta_0}
\subset Z_{r}^1\cup Z_{r}^2\cup Z_{\theta}^1\cup Z_{\theta}^2, \, \text{ with }
\\
&
Z_{r}^1:=\left\{ R-\mu/2\leq r \leq R-\nu/2,\ |\theta_{1}| \leq \nu \theta_0 \right\},\\
&
Z_{r}^2:=\left\{R+\nu/2\leq r \leq R+\mu/2,\ |\theta_{1}|\leq \nu \theta_0\right\},\\
&
Z_{\theta}^1:=\left\{R-\mu/2 \leq  r \leq R+\mu/2,\ -\mu \theta_0\leq\theta_{1}\leq-\nu \theta_0\right\},\\
&
Z_{\theta}^2:=\left\{R-\mu/2\leq r \leq  R+\mu/2,\ \nu \theta_{0}\leq\theta_{1}\leq \mu\theta_0\right\}.
\end{align*}

\noindent
\underline{On $Z_{r}^1.$} We use the structural hypothesis \eqref{conditionpetitpotentiel} to get
\begin{align}
\label{minoration1}
\nonumber
&
F_{r}(x)=-\lambda f'(r)g(\theta_{1})(r-R\cos(\theta_{1}))
\geq
-\lambda f'(r)g(\theta_{1})(R-\frac{\nu}{2}-R\cos(2\theta_{0}))
\\
&
\quad\,
\geq
\lambda f'(r)g(\theta_{1}) \frac{\nu-1}{2}
\geq 
\lambda (\nu-1)  \, \left(\min_{s \in  [-\mu,-\nu]} \chi'(s) \right)
\left(\min_{|s|\leq  \nu} \chi(s) \right)
=:c_1>0.
\end{align}
A similar proof establishes that, whenever $x \in Z_r^2$ we have
\begin{equation*}
F_{r}(x)\geq 
\lambda (\nu-1)  \, \left(\min_{s \in  [\nu,\mu]} [-\chi'(s)] \right)
\left(\min_{|s|\leq  \nu} \chi(s) \right)
=:c_2>0.
\end{equation*}

\noindent
\underline{On $Z_{\theta}^1$.} The important term is now $F_\theta$.
We have
\begin{align*}
&
F_{\theta}(x)
=
-\lambda\frac{R}{r}f(r)g'(\theta_{1})\sin(\theta_{1})
\geq
\lambda\frac{R}{r}f(r)g'(\theta_{1})\sin(\nu \theta_0)
\\
&
\quad\,
\geq
\lambda \frac{R}{\theta_0 (R+\mu/2)} \left( \min_{|s|\leq \mu} \chi(s) \right)
\left( \min_{s \in [-\mu,-\nu]} \chi'(s) \right)
=:c_3>0.
\end{align*}
A similar argument establishes that, whenever $x\in Z_\theta^2$ we have
\begin{align*}
F_{\theta}(x)
\geq
\lambda \frac{R}{\theta_0 (R+\mu/2)} \left( \min_{|s|\leq \mu} \chi(s) \right)
\left( \min_{s \in [\nu,\mu]} [-\chi'(s)] \right)
=:c_4>0.
\end{align*}

\medskip

Gathering all estimates,
there exists a positive constant $c_{\nabla}>0$
such that
\begin{equation*}
\forall\, x\,\in C_{R+\mu/2,\mu\theta_0}\setminus C_{R+\nu/2,\nu\theta_0},
\qquad
F_r(x)+F_\theta(x) \geq c_{\nabla}>0.
\end{equation*}

\paragraph{$\bullet$ Step five: end of the proof.}
Putting all  estimates together, we obtain
\begin{equation*}
\forall\, x\in\Pi_{x}H_{0}=\R^d\setminus B_\emptyset,\quad
F_r(x)+F_\theta(x) +n^2(x)\geq\min(c_{n^2},c_{\nabla})=:\alpha>0.
\end{equation*}
The lemma is proved.
\vspace{-0,4cm}
\hfill$\blacksquare$\\

\subsection{Refocusing Set}
\label{statset}

The goal of this subsection is to establish part (ii) of our main Theorem \ref{propnonconvergence}.

\medskip

Our main result is

\begin{prop}\label{lemstatset}Let $n^2$ be the potential defined in \eqref{n2pot}. Assume
the structural hypothesis \eqref{conditionpetitpotentiel} is fulfilled, namely $1-\cos(2\theta_0)<1/(2R)$.
Then, the refocusing set defined in Definition \ref{nonre} as
$$
M
=
\left\{(t,\xi,\eta)\in]0,+\infty[\times\R^{2d} \, \text{ s.t. } \,
\frac{|\eta|^2}{2}=n^2(0), \,
X(t,0,\xi)=0, \,
\Xi(t,0,\xi)=\eta\right\}
$$
satisfies
\begin{equation*}
M=\left\{(T_{R},\xi,\eta),\ \text{s.t.}\ \xi=-\eta=(r,\theta_{1},\ldots,\theta_{d-1}),\ r=\sqrt{2n^2(0)},\ |\theta_{1}|\leq\theta_{0}\right\},
\end{equation*}
where $T_{R}>0$ is the unique positive time such that
$X(T_{R},0,(\sqrt{2n^2(0)},0\ldots,0))=0$.
\end{prop}

\noindent
{\bf Proof of Proposition \ref{lemstatset}.}

Consider a trajectory $X(t,0,\xi)\equiv X(t)$ on the zero energy level,
with $\xi=(r,\theta_{1},\ldots,\theta_{d-1})$ in hyperspherical coordinates.

If $|\theta_1|\geq 2\theta_{0}$, it is clear that $X(t)$ is a straight line which never enters  $B_{p}$, and
the equation $X(t,0,\xi)=0$ with $t>0$ has no solution.

We need to understand the geometry when the trajectory reaches $B_{p}$, {\em i.e.}
when $|\theta_1|< 2\theta_0$.

We prove below that two cases occur. If $|\theta_{1}|\leq\theta_{0}$, the trajectory remains along a line, and it is reflected by
the refraction index towards the origin.
If $\theta_0<|\theta_1|< 2 \theta_0$, the force acting on the trajectory has a non-vanishing component in
the orthoradial direction, which prevents the trajectory to go back to the origin. The proposition  follows.

Let us come to a proof.

\medskip

$\bullet$ \textbf{First case: $|\theta_1|\leq \theta_0$.}

Consider the trajectory $Y(t)$ defined in hyperspherical coordinates as
\begin{equation*}
Y(t)=(r(t),\theta_{1},\ldots,\theta_{d-1})\,,
\end{equation*}
with $r(t)$ solution to the  ordinary equation $r''=-\lambda f'(r)$ with initial data
\begin{equation*}
r(0)=0,\qquad r'(0)=\sqrt{2n^2(0)}.
\end{equation*}
Then, $(Y(t),Y'(t))$  satisfies the Hamiltonian ODE (\ref{eqHJ}) associated with $h(x,\xi)=|\xi|^2/2+n^2(x)$.
Since $Y(0)=X(0)=0$, and $ Y'(0)= X'(0)=\xi$, uniqueness provides $X(t)=Y(t)$ for all $t$.
The trajectory $X(t)$ is radial. 

It is clear that the radial trajectory $t\mapsto r(t)$ reaches the  region $\{ R-1\leq r\leq R+1\}$ at
time $t_{e}=(R-1)/|\xi|=(R-1)/\sqrt{2 \, n^2(0)}>0$, where
$
t_{e}=\inf\left\{t>0, \, X(t)\in B_{p}\right\}.
$
Now, according to Corollary \ref{cornontrapping}, the trajectory $r(t)$ necessarily leaves the 
region $\{ R-1\leq r\leq R+1\}$ at some later time $t_s> t_e$,
where
$
t_{s}=\inf\left\{t>t_e, \, X(t)\notin B_{p}\right\}.
$
The trajectory can either leave the bump at $r=R-1$ or at $r=R+1$. The case $r=R+1$ is forbidden, for in the contrary case,
using continuity, there would exist a time $t_{c}$ such that $r(t_{c})=R$, hence $X(t_{c})\in B_{\emptyset}$, which is not allowed.
Therefore, the trajectory leaves the bump $B_{p}$ at $X(t_s)$ where
 $|X(t_s)|=r(t_{s})=R-1$. 
Energy conservation, together with  the fact that the trajectory is radial, implies that $X'(t_s)=-\xi$.
Therefore, 
the trajectory for later times $t\geq t_s$ is  a straight line with constant speed $-\xi$. We deduce
that there exists a unique $T_R>t_s$ such that $X(T_R,0,\xi)=0$, and we have as desired
$\Xi(T_R,0,\xi)=-\xi$.

\medskip

$\bullet$ \textbf{Second case: $\theta_0<|\theta_1|< 2\theta_0$.}

We first assume that $d=2$, and next generalize the argument to $d\geq 3$ using the symmetries of the system.
To fix the ideas, we assume in the following that $\theta_0<\theta_1< 2\theta_0$, the proof being the same when $\theta_1$ has the opposite sign.

\mbox{}\qquad
{\bf $*$ In dimension $d=2$.}

Let $t_{e}=(R-1)/|\xi|$ be the time when the trajectory enters $B_{p}$, as in the preceding case.

On the one hand, since the velocity $\Xi(t_e)$ is radial and satisfies
$\Xi(t_e)=|\xi| \, \vec{u}_{r}$, there is an $\eps>0$ such that $R-1<|X(t)|<R+1$ whenever $t\in ]t_e,t_e+\eps]$.
On the other hand, by assumption we have $\theta_1(t_{e})=\theta_1\in]\theta_{0},2\theta_0[$, and continuity
implies there is an $\eps>0$ such that $\theta_0<\theta_1(t)<2 \theta_0$
whenever $t\in [t_e,t_e+\eps]$.
Hence  we may  define
$$
t_{s}:=\sup\{t \geq t_{e}, \,  \, \text{ s.t. } \, \forall t' \in [t_e,t], \quad \theta_1(t')\in ]\theta_0,2\theta_0[ \text{ \, { and } \, }
X(t')\neq 0.\}.
$$
Now, Hamilton's equations of motion  \eqref{eqHJ} can be written in polar coordinates as
\begin{equation*}
\begin{cases}
r''-r(\theta_{1}')^2=-\lambda f'(r)g(\theta_{1}),\\
2r'\theta_{1}'+r\theta_{1}''=-\lambda \frac{f(r)}{r}g'(\theta_{1}).
\end{cases}
\end{equation*}
Examining the second equation, we have $(r^2\theta_{1}')'=2rr'\theta_{1}'+r^2\theta_{1}''=-\lambda f(r) g(\theta_{1})$, 
and we get whenever $r(t)\neq0$, 
\begin{equation}
\label{eqmachintruc}
\theta_{1}'(t)=-\frac{\lambda}{r^2(t)}\int_{t_{e}}^tf(r(s))g'(\theta_{1}(s))ds.
\end{equation}
Therefore, since  $f(r)\geq 0$ for any $r\geq 0$ while $f(r)>0$ whenever $R-1<r<R+1$,
and  since  $g'(\theta_1)\leq 0$ when $\theta_0\leq \theta_1 \leq 2\theta_0$,
while $g'(\theta_1)< 0$ when $\theta_0< \theta_1 < 2\theta_0$
we get, with the above definitions and observations,
$$
\theta_1'(t)>0 \quad \forall t\in]t_e,t_s].
$$
With this observation at hand, two cases may occur.

If $t_s=+\infty$, there is nothing to prove, for by definition of $t_s$, we have $X(t)\neq 0$ whenever
$0<t\leq t_s=+\infty$.

In the case $t_s<+\infty$, we already know $X(t)\neq 0$ whenever $0<t\leq t_s$. Besides, since
$\theta_1'(t)>0$ whenever  $0<t\leq t_s$, it is clear that the case $X(t_s)=0$ is impossible 
(for in that case the trajectory would be a straight line passing through the origin on some interval $[t_*,t_s]$,
in contradiction with $\theta_1'(t)>0$ on  $[t_*,t_s]$), hence
$\theta_1(t_s)=2\theta_0$ and  $\theta_1'(t_s)>0$.
For that reason, the trajectory $X(t)$ for times $t> t_s$ is a straight line with constant velocity, which 
lies entirely in the set $2\theta_0<\theta_1<2\theta_0+\pi$. In particular, since $\theta_1'(t_s)>0$, the trajectory cannot be radial and we have $X(t)\neq 0$ whenever $t>t_s$ in that case. This concludes the proof.
\medskip

\mbox{}\qquad
{\bf $*$ In dimension $d\geq 3$}.

We use the invariance of $n^2$ under the action of $\Ot_{d,1}(\R)$.

Take $\xi\in\R^d$ such that $|\xi|=\sqrt{2n^2(0)}$.
Write $\xi=(\sqrt{2n^2(0)},\theta_{1},\ldots,\theta_{d-1})$ in hyperspherical coordinates.
There exists a matrix $A_\xi \in \Ot_{d,1}(\R)$ such that $A_\xi \, \xi=(\sqrt{2n^2(0)},\theta_{1},0,\ldots,0)$. 
On the other hand,
denote by $(r(t),\theta_{1}(t))$ the solution of  Hamilton's equations of motion
\eqref{eqHJ} with initial data $(\sqrt{2n^2(0)},\theta_{1})$ in dimension $2$. 
 We set $Y(t)=A_\xi^{-1}(r(t),\theta_{1}(t),0\ldots,0)$. Then $Y(t)$ satisfies Hamilton's equations of motion \eqref{eqHJ}, with initial data $Y(0)=0,\ Y'(0)=\xi$. Uniqueness  provides $Y(t)=X(t)$ for any $t> 0$.
This, combined with the previous step, provides $X(t)\neq 0$ for any $t> 0$.
 \vspace{-0,4cm}
\hfill $\blacksquare$

\bigskip


\section{Convergence proof}
\label{proposition}


The goal of this section is to prove item (iii) of our main Theorem \ref{propnonconvergence}.

The proof is performed in a number of steps.
We begin by defining some necessary  notation.


\subsection{The linearized hamiltonian flow}


Let $\varphi(t,x,\xi)=\left(X(t,x,\xi),\Xi(t,x,\xi)\right)$ denote the flow associated with Hamilton's equations of motion
\eqref{eqHJ}.
The linearized flow, written $F(t,x,\xi)$ below, is
\begin{equation*}
F(t,x,\xi)= \frac{D\varphi(t,x,\xi)}{D(x,\xi)}:=\begin{pmatrix}A(t,x,\xi) & B(t,x,\xi)\\ C(t,x,\xi)&D(t,x,\xi)\end{pmatrix},
\end{equation*}
where $A(t)$, $B(t)$, $C(t)$, $D(t)$ are by definition
\begin{align*}
&A(t,x,\xi)=\frac{DX(t,x,\xi)}{Dx},\ \ B(t,x,\xi)=\frac{DX(t,x,\xi)}{D\xi},\\
&C(t,x,\xi)=\frac{D\Xi(t,x,\xi)}{Dx},\ \ D(t,x,\xi)=\frac{D\Xi(t,x,\xi)}{D\xi}.
\end{align*}
The linearisation of $(\ref{eqHJ})$ leads to
\begin{equation}
\label{eqHJL1}
\left\{
\begin{array}{ll}
\displaystyle \frac{\partial}{\partial t}A(t,x,\xi)=C(t,x,\xi),& A(0,x,\xi)=Id,\\
\displaystyle \frac{\partial}{\partial t}C(t,x,\xi)=\frac{D^2n^2}{Dx^2}(X(t,x,\xi))A(t,x,\xi),&C(0,x,\xi)=0,
\end{array}
\right.
\end{equation}
 as well as
\begin{equation}
\label{eqHJL2}
\left\{
\begin{array}{ll}
\displaystyle
\frac{\partial}{\partial t}B(t,x,\xi)=D(t,x,\xi),& B(0,x,\xi)=Id,\\
\displaystyle
 \frac{\partial}{\partial t}D(t,x,\xi)=\frac{D^2n^2}{Dx^2}(X(t,x,\xi))B(t,x,\xi),&D(0,x,\xi)=0.
\end{array}
\right.
\end{equation}
Finally, we define for later purposes the matrix $\Gamma(t,x,\xi)$ as
\begin{equation}
\label{gamma}
\Gamma(t,x,\xi)=\left( C(t,x,\xi)+i \, D(t,x,\xi)\right) \, . \, \left( A(t,x,\xi)+i B(t,x,\xi)\right)^{-1}.
\end{equation}


\subsection{A wave packet approach: preparing for a stationary phase  argument}


The intermediate result in Proposition \ref{intermres} establishes roughly that
$\langle w_\eps,\phi\rangle \sim \langle w_{out}+\widetilde{w_\eps},\phi\rangle$
as $\eps\to 0$.  Therefore, item (iii) of our main Theorem reduces to proving
$\langle\widetilde{w_\eps},\phi\rangle \sim \langle L_\eps,\phi\rangle$ as $\eps \to 0$.

Therefore, this preliminary paragraph is devoted to
express the quantity
\begin{equation*}
\langle\widetilde{w_\eps},\phi\rangle
=\frac{1}{\varepsilon}\int_{\theta}^{T_1}\left(1-\chi\left(\frac{t}{\theta}\right)\right)e^{-\alpha_\varepsilon t}\left\langle U_\varepsilon(t)\chi_\delta(H_\varepsilon)S_\varepsilon,\phi_\varepsilon\right\rangle dt.
\end{equation*}
as an appropriate oscillatory integral.
Our approach uses the technique developped in 
\cite{MR2139886}, which in turn strongly uses a wave packet theorem
due to M. Combescure and D. Robert (see \cite{MR1461126}). We skip here the details of the proof, referring to
\cite{MR2139886}.

The main result in this paragraph is the following
\begin{prop}
\label{reduction}
{\bf (See \cite{MR1461126})} \,
Whenever $X=(q,p,x,\xi,y,\eta)\in \R^{6d}$ and $t\in \R$, define the complex phase
\begin{align}
&
\label{phase}
\nonumber
\psi(t,X):=\int_{0}^t\left(\frac{p_{s}^2}{2}+n^2(q_{s})\right)ds-p.(x-q)+p_{t}.(y-q_{t})
\\
&
\qquad\qquad\qquad
+x.\xi-y.\eta+i\frac{(x-q)^2}{2}+\frac{1}{2}\Gamma_{t}(y-q_{t}).(y-q_{t}),
\end{align}
where $q_t:=X(t,q,p)$, $p_t:=\Xi(t,q,p)$, and $\Gamma_t:=\Gamma(t,q,p)$.
Select an integer $N\in\N$.
Select two truncation functions $\chi_0(q,p)$ and $\chi_1(x,y)$ both lying in $C_0^\infty(\R^{2d})$, and
such that
\begin{align*}
&
\text{supp}\ \chi_{0}(q,p)\subset\left\{|q|\leq2\delta\right\}\cup\left\{||p|^2/2-n^2(q)|\leq 2\delta\right\}\,,
\\
&
\chi_{0}(q,p)\equiv 1\ \text{on}\ \left\{|q|\leq3\delta/2\right\}\cup\left\{||p|^2/2-n^2(q)|\leq 3\delta/2\right\},
\\
&
\chi_{1}(x,y)\equiv 1 \, \text{ close to } \, (0,0).
\end{align*}
Define the amplitude
\begin{align}
\label{amplitude}
a_{N}(t,X):=
e^{-\alpha_{\varepsilon}t}(1-\chi)\left(\frac{t}{\theta}\right)
\widehat{S}(\xi)\widehat{\phi}^*(\eta)\chi_{0}(q,p)\chi_{1}(x,y)P_{N}\left(t,q,p,\frac{y-q_{t}}{\sqrt\varepsilon}\right),
\end{align}
where $P_N(t,q,p,z)$ satisfies
\begin{align}
\label{pn}
P_{N}(t,q,p,x):=\frac{1}{\pi^{d/4}}det(A(t,q,p)+iB(t,q,p))_c^{-1/2}\mathcal Q_{N}(t,q,p,x),
\end{align}
and the square root $det(A(t,q,p)+iB(t,q,p))_c^{-1/2}$ is defined by continuously following the argument of
the relevant complex number, starting from the value $det(A(0,q,p)+iB(0,q,p)=1$ at time $t=0$, while
$\mathcal Q_{N}(t,q,p,x)$ is a polynomial in the variable $x\in \R^d$, whose coefficients vary smoothly with
$(t,q,p)$, and $\eps$, and which satisfies
$$
\mathcal Q_{N}(t,q,p,x)=1+O(\sqrt{\eps})
$$
in the relevant topology. More precisely, we have
\begin{equation}\label{eqQN}
\begin{cases}
\mathcal Q_N(t,q,p,x)=1+\displaystyle\sum_{(k,j)\in I_N}\varepsilon^{\frac{k}{2}-j}p_{k,j}(t,q,p,x),\\
I_N=\left\{1\leq j\leq 2N-1,1\leq k-2j\leq 2N-1, \, k\geq 3j\right\},
\end{cases}
\end{equation}
where each $p_{k,j}$ has at most degree $k$ in the variable $x$.
\medskip

Then, the following holds
\begin{align}
\label{oscinteg}
\langle\widetilde{w_\eps},\phi\rangle
=
\frac{1}{\varepsilon^{(5d+2)/2}}\int_{\theta}^{T_{1}}\int_{\R^{6d}}e^{\frac{i}{\varepsilon}\psi(t,X)}a_{N}(t,X)dtdX +O_{T_{1},\delta}(\varepsilon^N).
\end{align}
\end{prop}

\medskip

\noindent
{\bf Sketch of proof of Proposition \ref{reduction}.}

\noindent
Using the short-hand notation
$\widetilde\chi_{\delta}(t):=e^{-\alpha_{\varepsilon}t}(1-\chi)\left(t/\theta\right)$,
we have
\begin{align}
\label{wep}
\langle\widetilde{w_\eps},\phi\rangle
=i/\varepsilon\int_{\theta}^{T_{1}}
\widetilde\chi_{\delta}(t) \,
\left\langle \chi_{\delta}(H_{\varepsilon})S_{\varepsilon},U_{\varepsilon}(-t)\phi_{\varepsilon}\right\rangle dt.
\end{align}
To compute the term
$U_{\varepsilon}(-t)\phi_{\varepsilon}$ accurately, we use a projection over the overcomplete basis
of~$L^2(\R^d)$  obtained by using the so-called gaussian wave-packets, namely the family of functions indexed by $(q,p)\in\R^{2d}$ defined by
\begin{equation*}
\varphi_{q,p}^\varepsilon(x,\xi):=\frac{1}{(\pi\varepsilon)^{d/4}}\exp\left(\frac{i}{\varepsilon}p.\left(x-\frac q 2\right)\right)\exp\left(-\frac{(x-q)^2}{2\varepsilon}\right).
\end{equation*}
The point indeed is that, as proved by Combescure and Robert in
\cite{MR1461126}, we have
\begin{align}
\label{combrob}
\nonumber
&
U_{\varepsilon}(-t)\,
\varphi_{q,p}^\varepsilon(x,\xi)
=O_{T_1,\delta}(\eps^N)+
\\
&
\nonumber
\qquad
\frac{1}{\varepsilon^{d/4}}
\exp\left(\frac{i}{\varepsilon}p_t.\left(x-\frac{q_t}{2}\right)\right) \,
\exp\left(  -\frac{|x-q_t|^2}{2\eps}\right) \,
\\
&
\qquad\quad
\exp\left(\frac{i}{\eps}
\left[
\int_0^t \left(\frac{p_s^2}{2}+n^2(q_s)\right) ds
-\frac{q_t \cdot p_t - q \cdot p}{2}
\right]
\right)
\,
P_N\left(t,q,p,\frac{x-q_t}{\sqrt{\eps}}\right) \,
\end{align}
in $L^\infty([0,T_1];L^2(\R^{d}))$. In other words, we have a quite explicit complex-phase/amplitude representation
of the Schr\"odinger propagator when acting on the gaussian wave packets.

This observation leads to writing, successively, in (\ref{wep})
\begin{align*}\
&
\left\langle \chi_{\delta}(H_{\varepsilon})S_{\varepsilon},U_{\varepsilon}(-t)\phi_{\varepsilon}\right\rangle
=
\frac{1}{(2\pi\varepsilon)^d}\int_{\R^{2d}}\left\langle\chi_{\delta}(H_{\varepsilon})S_{\varepsilon},\varphi_{q,p}^\varepsilon\right\rangle\left\langle\varphi_{q,p}^\varepsilon,U_{\varepsilon}(-t)\phi_{\varepsilon}\right\rangle dq\, dp\,,
\\
&
\qquad
\qquad
\qquad
\qquad
\quad\,\,\,\,
=
\frac{1}{(2\pi\varepsilon)^d}\int_{\R^{2d}}
\left\langle\chi_{\delta}(H_{\varepsilon})S_{\varepsilon},\varphi_{q,p}^\varepsilon\right\rangle
\,
\left\langle U_{\varepsilon}(t)\varphi_{q,p}^\varepsilon,\phi_{\varepsilon}\right\rangle dq\, dp\,.
\end{align*}
Now, the idea is to replace the factor $U_{\varepsilon}(t)\varphi_{q,p}^\varepsilon$ by its approximation derived above.
Yet a few preliminary steps are in order. The first one
uses the truncation in energy $\chi_{\delta}(H_{\varepsilon})$, together with the functional calculus for pseudo-differential operators of Helffer and Robert
(see \cite{MR724029} ), to replace this truncation by an explicit truncation near the set $p^2/2+n^2(q)=0$, up to small
error terms.
The second step consists in using the Parseval formula to write (we want to exploit the source term $S_\eps$ on the Fourier side)
\begin{align*}
\langle S_\eps \, , \phi_{q,p}^\eps \rangle
=
\frac{1}{(2\pi\eps)^{d/2}} \int e^{i\frac{x \cdot \xi}{\eps}} \, \widehat{S}(\xi) \, \phi_{q,p}^\eps(x) \, dx \, d\xi
=
\frac{1}{(2\pi\eps)^{d/2}} \int \widetilde\chi(x) \,e^{i\frac{x \cdot \xi}{\eps}} \, \widehat{S}(\xi) \, \phi_{q,p}^\eps(x) \, dx \, d\xi,
\end{align*}
for some function $\widetilde\chi(x)$ that truncates close to $x=0$, and similarly
\begin{align*}
\left\langle U_{\varepsilon}(t)\varphi_{q,p}^\varepsilon,\phi_{\varepsilon}\right\rangle
=
\frac{1}{(2\pi\eps)^{d/2}} \int \widetilde\chi(y) \,e^{i\frac{y \cdot \eta}{\eps}} \, \widehat{\phi}(\eta) \,
\left( U_{\varepsilon}(t) \varphi_{q,p}^\eps\right) (y) \, dy \, d\eta,
\end{align*}
These two steps explain the truncation factors $\chi_0$ and $\chi_1$ in the Proposition,
which act close to the zero energy-level in phase-space (this is where functional calculus is used) and close to the origin 
in physical space. The last step consists in exploiting formula (\ref{combrob}) in the obtained representation.

Eventually, one obtains the desired formula.
 \vspace{-0,4cm}
\hfill $\blacksquare\\$


\subsection{Preparing for a stationary phase argument}


This slightly technical paragraph is devoted to proving that the obtained phase $\psi$ in Proposition~\ref{reduction} satisfies the assumptions of the stationary phase Theorem.

Our main result in this paragraph is the Proposition after the following Lemma.

\begin{lem}
\label{propM}
Let $n^2$ be any smooth refraction index.
Then, the following holds

\medskip

(i) The stationary set associated with the phase $\psi$ in (\ref{phase}), defined as
\begin{align*}
M_X:=
\left\{
(t,X)=(t,q,p,x,\xi,y,\eta)\in[\theta,T_1]\times\R^{6d} \,  \text{ s.t. } \,
\nabla_{t,X}\psi(t,X)=0 \, \text{ and } \, {\rm Im} \, \psi(t,X)=0
\right\}
\end{align*}
satisfies
\begin{align}
\label{stM}
M_X
=
\{ (t,q,p,x,\xi,y,\eta) \, \text{ s.t. } \, x=y=q=0, \, \xi=p, \, (t,p,\eta) \in M\},
\end{align}
where we recall that $M=\{(t,p,\eta), \, X(t,0,p)=0, \, \Xi(t,0,p)=\eta, \, \eta^2/2=n^2(0)\}$ by definition.

\medskip

(ii)
We have, whenever $m=(t,X)\in M_X$, the relation
\begin{align}
\label{ker}
&{\rm Ker}(D^2\psi_{|_{m}})=
\left\{(T,Q,P,X,\Xi,Y,H)\in]0,+\infty[\times\R^{6d},\ X=Y=Q=0,\right.
\\
\nonumber
&\hspace{3cm}\label{ker} \left. \Xi=P,\ \eta^TH=0,B_t(0,p) \, P+T\eta=0,\ -H+D_t(0,p) \,  P+T\nabla n^2(0)=0\right\}.
\end{align}
\end{lem}

\medskip

Note that this Lemma does not use the particular structure of our index.

\medskip

\noindent
{\bf Proof of Lemma \ref{propM}.}
A mere computation of $\text{Im}\, \psi$ and $\nabla\,\psi$ allows to write $(\ref{stM})$. Differentiating $\nabla\,\psi$ once allows to write $(\ref{ker})$. For more details, the reader may check $\cite{MR2139886}$.
\vspace{-0,4cm}
\hfill $\blacksquare$\\

\medskip

With this Lemma at hand, our key result in this section is the following 
 
\medskip

\begin{prop}
\label{h2}
Let $n^2$ be the refraction index defined in \eqref{n2pot}.
We recall that the refocusing set $M$ is computed in Lemma \ref{lemstatset}
and satisfies
$$
M=\left\{(T_R,\xi,\eta) \, \text{ s.t. } \,
\xi=-\eta=(r,\theta_1,\ldots,\theta_{d-1}), \, r=\sqrt{2 n^2(0)}, \, |\theta_1|\leq \theta_0\right\}.
$$
Now, take any
\begin{align*}
&
m \in \overset{\circ}{M}_{X}=\Big\{(t,q,p,x,\xi,y,\eta) \, \text{ s.t. } \, x=y=q=0, \, \xi=p, \,
\\
&
\qquad\qquad\qquad
(t,p,\eta) \in M, \, \text{with} \, p=(r,\theta_1,\ldots,\theta_{d-1}), \, \text{ and }  \,  |\theta_1|<\theta_0\Big\}
\end{align*}
Then, we have
$$
{\rm Ker} \, D^2\psi|_m=T_m M_X,
$$
where $T_m M_X$ denotes the  space tangent to $M_x$ at point $m$.
\end{prop}

The remainder part of this subsection is devoted to the proof of Proposition \ref{h2}.
We begin by proving the Proposition in the
case
$$
m=m_{0}:=(T_R,0,p_0,0,p_0,0,-p_{0}),\quad \text{where}\quad p_{0}:=(\sqrt{2n^2(0)},0,\ldots,0)\,.
$$ 
We next generalize the result to other values of $m$, using the symmetries of the problem.


\subsubsection{Proof of Proposition \ref{h2} when $m=m_{0}$}


The computation of $T_{m_{0}}M_{X}$ on the one hand is rather easy
\begin{lem} The space $T_{m_0} M_{X}$ is given by
\label{lemtang}
\begin{equation*}
T_{m_{0}}M_{X}=\{(T,Q,P,X,\Xi,Y,H)\ \text{s.t.}\ X=Y=Q=T=0,\ \Xi=P=-H,\ P.p_{0}=0\}.
\end{equation*} 
\end{lem}

\medskip

\noindent
{\bf Proof of Lemma \ref{lemtang}.}
This is a mere computation starting from the definition of the refocusing set $M$, as
$M=\{(t,p,\eta), \, X(t,0,p)=0, \, \Xi(t,0,p)=\eta, \, \eta^2/2=n^2(0)\}$.
\vspace{-0,3cm}
\hfill$\blacksquare$\\

\bigskip

In order to determine ${\rm Ker} \, D^2\psi_{|_{m_{0}}}$ the first step it to compute the matrices $B_t$ and $D_t$
involved in the linearized flow, see \eqref{eqHJL2}.
\begin{lem}
\label{lemlin}
Let $n^2$ be the potential defined in \eqref{n2pot}. Then, we have
\begin{equation}\label{eqlinBD}\hspace{-0,29cm}
D(T_{R},0,p_{0}):=\frac{\partial\Xi}{\partial\xi}(T_{R},0,p_{0})=-I_{d}, \quad B(T_{R},0,p_{0}):=\frac{\partial X}{\partial\xi}(T_{R},0,p_{0})=\begin{pmatrix}b_{11}&0&\\0&O_{d-1}\end{pmatrix},
\end{equation}
where $I_{d}$ is the identity matrix, $b_{11}\in\R$ and $O_{d-1}$ is a square matrix of dimension $d-1$ equal to~$0$.  

\end{lem}

\noindent
{\bf Proof of Lemma \ref{lemlin}.}

 We consider $x_{0}(t,0,p)=(x_{0}^1(t,0,p),\ldots,x_{0}^d(t,0,p))$ the solution to~\eqref{eqHJ} with initial data $x_{0}(0,0,p)=0$ and $x'_0(0,0,p)=p$. 

We recall that the index $n^2$ is invariant under the action of $\Ot_{d,1}(\R^d)$. Thus we first compute the components
of $D$ and $B$ that are
invariant under $\Ot_{d,1}(\R^d)$, namely their first column. We next compute the other columns
by using the symmetries again, in conjunction with a perturbation argument.

\medskip

$\bullet$ \textbf{Computation of $\frac{\partial \Xi}{\partial\xi_{1}}(T_{R},0,p_{0})$ and $\frac{\partial X}{\partial\xi_{1}}(T_{R},0,p_{0})$}

 We start with $\frac{\partial \Xi_{j}}{\partial\xi_{1}}(T_{R},0,p_{0})$ for $j\geq 2$. We have
\begin{align*}
\frac{\partial \Xi_{j}}{\partial\xi_{1}}(T_{R},0,p_{0})
&=
\lim_{\varepsilon\rightarrow0}
\frac{\Xi_{j}\left(T_{R},0,(\sqrt{2n^2(0)}+\varepsilon,0\ldots,0)\right)
-
\Xi_{j}\left(T_{R},0,(\sqrt{2n^2(0)},0\ldots,0)\right)}{\varepsilon}.
\end{align*}
Since the trajectory is radial we have 
 \begin{equation*}
\Xi_{j}\left(T_{R},0,(\sqrt{2n^2(0)}+\varepsilon,0\ldots,0)\right)=\Xi_{j}\left(T_{R},0,(\sqrt{2n^2(0)},0\ldots,0)\right)=0,\quad \forall\ j\geq2.
\end{equation*}
Hence,
$\frac{\partial \Xi_{j}}{\partial\xi_{1}}(T_{R},0,p_{0})=0$, $\forall\ j\geq2$. A similar argument provides $\frac{\partial X_{j}}{\partial\xi_{1}}(T_{R},0,p_{0})=0$, $\forall\ j\geq2$. There remains to determine the first coefficient of $D$,
namely $\frac{\partial \Xi_{1}}{\partial\xi_{1}}(T_{R},0,p_{0})$.
Since the trajectory is radial, and by conservation of the energy, we have for $\varepsilon$ small enough
\begin{align*}
&
\Xi_{1}\left(T_{R},0,(\sqrt{2n^2(0)}+\varepsilon,0,\ldots,0)\right)
=-\left(\sqrt{2n^2(0)}+\varepsilon\right),
\\
&
\Xi_{1}\left(T_{R},0,(\sqrt{2n^2(0),0,\ldots,0)}\right)=-\sqrt{2n^2(0)}.
\end{align*}
Thus,
\begin{equation*}
d_{11}:=\lim_{\varepsilon\rightarrow 0^+}\frac{\Xi\left(T_{R},0,(\sqrt{2n^2(0)}+\varepsilon,0,\ldots,0)\right)
-\Xi\left(T_{R},0,(\sqrt{2n^2(0)},0,\ldots,0)\right)}{\varepsilon}=-1.
\end{equation*}

\medskip

$\bullet$ \textbf{Computation of $\frac{\partial \Xi(T_{R},0,p_{0})}{\partial\xi_{j}}$ and $\frac{\partial X(T_{R},0,p_{0})}{\partial\xi_{j}}$ ($j\geq2$)}

 Considering the symmetries of the problem, it is enough to
consider the case $j=2$: the other components may be determined using the same argument.

\begin{figure}[htb]
\begin{center}
\includegraphics[scale=0.4]{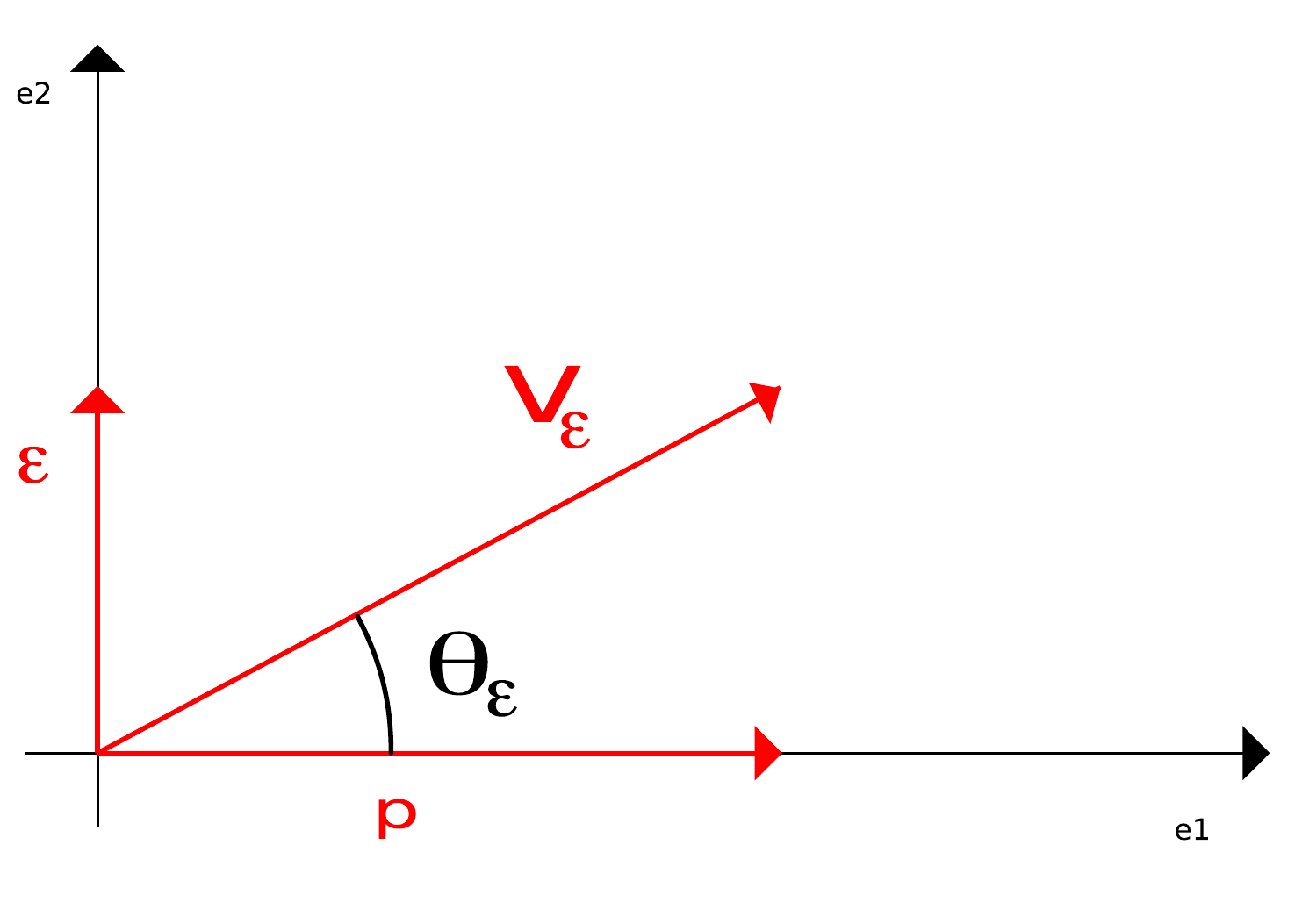}
\caption{Perturbation of the initial speed}
\label{Perturbation}
\end{center}
\end{figure}

We perturb the initial speed along the direction $e_{2}$, by a factor $\eps$ (see Figure $\ref{Perturbation}$).

Let $X_{\varepsilon}(t)$ be the solution of the perturbed problem
\begin{equation*}
\begin{cases}
X_{\varepsilon}''(t)=\nabla n^2(X_{\varepsilon}(t)),\qquad
X_{\varepsilon}(0)=0,\qquad X_{\varepsilon}'(0)=p_{0}+\varepsilon \,e_{2}.
\end{cases}
\end{equation*}
We expand $X_{\varepsilon}(t)$ with respect to $\varepsilon$ and obtain $X_{\varepsilon}(t)=X_{0}(t)+\varepsilon X_{1}(t)+\ldots$. With this notation we have $X_{1}(t)=\frac{\partial X}{\partial\xi_{2}}(t)$
and $X_{1}'(t)=\frac{\partial\Xi}{\partial\xi_{2}}(t)$.
To obtain the expansion in $\varepsilon$, we go back to the previous case ($j=1$) using a change of variables. Indeed,
for $\varepsilon$ small enough, the trajectory is radial along the direction $X_{\varepsilon}'(0)$. Let $(\widetilde{e_{1}}, \ldots, \widetilde{e_{d}})$ be a new basis defined by $\widetilde{e_{j}}:=O_{\varepsilon}e_{j}$,
with 
\begin{equation*}
O_{\varepsilon}:=
\begin{pmatrix}
\cos(\theta_{\varepsilon})& -\sin(\theta_{\varepsilon})&0&\ldots&0\\
\sin(\theta_{\varepsilon})&\cos(\theta_{\varepsilon})&0&\ldots&0\\
0&0\\
\vdots&\vdots&&I_{d-2}\\
0&0
\end{pmatrix},\qquad \cos(\theta_{\varepsilon})=\frac{p_{0}}{p_{0}^2+\varepsilon^2},\qquad \sin(\theta_{\varepsilon})=\frac{\varepsilon}{p_{0}^2+\varepsilon^2}.
\end{equation*}
 Let $\widetilde{X_{\varepsilon}}$ be the coordinates of $X_{\varepsilon}$ in $(\widetilde{e_{1}},\ldots,\widetilde{e_{d}})$. Since $O_{\varepsilon}^{-1}\nabla n^2(X_{\varepsilon})=\nabla n^2(\widetilde{X_{\varepsilon}})$, we clearly have
\begin{equation*}
\widetilde{X_{\varepsilon}}''(t)=\nabla n^2(\widetilde{X_{\varepsilon}}(t)),\quad\widetilde{X_{\varepsilon}}(0)=0,\quad \widetilde{X_{\varepsilon}}'(0)=(\sqrt{\varepsilon^2+p_{0}^2},0,\ldots,0)=p_{0}+O(\varepsilon^2).
\end{equation*}
Hence it is clear that $\widetilde{X_{\varepsilon}}(t)=\widetilde{X_{0}}(t)+O(\varepsilon^2)$.
Therefore, we recover
\begin{align*}
X_{0}(t)+\varepsilon X_{1}(t)&=O_{\varepsilon}\left(\widetilde{X_{0}}(t)+O(\varepsilon^2)\right)=(I_{d}+\varepsilon E+O(\varepsilon^2))(\widetilde{X_{0}}(t)+O(\varepsilon^2)),
\end{align*}
with 
\begin{equation*}
E:=\begin{pmatrix}
0& -\frac{1}{p_{0}}&0&\ldots&0\\
\frac{1}{p_{0}}&0&0&\ldots&0\\
0&0\\
\vdots&\vdots&&I_{d-2}\\
0&0
\end{pmatrix}.
\end{equation*}
In other words, we have
\begin{equation*}
\forall\, t\in\R,\qquad X_{0}(t)=\widetilde{X_{0}}(t)\qquad \text{and}\qquad X_{1}(t)=E \widetilde{X_{0}}(t).
\end{equation*}
Since the Hamiltonian trajectory goes back to the origin at time $T_{R}$, we deduce
\begin{equation*}
\frac{\partial X}{\partial\xi_{2}}(T_{R},0,p_{0})=X_{1}(T_{R})=E\widetilde{X_{0}}(T_{R},0,p_{0})=E\times0=0.
\end{equation*}
In the same way, we have
\begin{align*}
\frac{\partial\Xi}{\partial\xi_{2}}(T_{R},0,p_{0})&=X_{1}'(T_{R})=E\widetilde{X_{0}}'(T_{R})={}^t\left(-\frac{x_{0,R}'^2(T_{R})}{p_{0}},\frac{x_{0,R}'^1(T_{R})}{p_{0}},x_{0,R}'^3(T_{R}),\ldots,x_{0,R}'^d(T_{R})\right)\,,\\
&={}^t(0,1,0,\ldots,0).
\end{align*}
The columns of $B$ and $D$ (for $j\geq3$) are determined in the similar way.
This leads to \eqref{eqlinBD}.
\vspace{-0,3cm}
\hfill$\blacksquare$\\

\bigskip

At this stage, we deduce the
 \begin{cor}
 \label{proptanker0}
 $\text{Ker}\, D^2\psi_{|_{m_{0}}}=T_{m_{0}}M_{X}.$
 \end{cor}
 
\noindent
{\bf Proof of Corollary \ref{proptanker0}.}

According to $(\ref{ker})$, we have
 \begin{multline*}
 Ker(D^2\psi_{|_{m}})=
\left\{(T,Q,P,X,\Xi,Y,H), \, X=Y=Q=0,\right.\\
\left.\Xi=P,\eta^TH=0,B_{T_{R}}(0,p) \, P+T\eta=0,-H+D_{T_{R}}(0,p) \, P+T\nabla n^2(0)=0\right\}.
\end{multline*}
 Since $\eta=-p_{0}$, we recover $H=(0,H_{2},\ldots,H_{d})$ (in Cartesian coordinates). Since $\nabla n^2(0)=0$, we deduce that $D_{T_{R}}(0,p) \, P=H$.
According to Lemma $\ref{lemlin}$, we deduce that $H=-P$. Finally, $B_{T_{R}(0,p) \, }P=0$ hence $T=0$. Thus,
 \begin{equation*}
  \text{Ker}\, D^2\psi_{|_{m_{0}}}=
\left\{(T,Q,P,X,\Xi,Y,H), \, X=Y=Q=T=0,\ P=\Xi=-H,\ P.p_{0}=0\right\}.
 \end{equation*}
Using Lemma \ref{lemtang}, the proof is complete.
 \vspace{-0,3cm}
 \hfill $\blacksquare$\\


\subsubsection{Proof of Proposition  \ref{h2} for any $m$}
 

 In this subsection, we prove the
 \begin{lem}
 \label{lemH2}
$\forall\, m\in \overset{\circ}{M}_{X}$, we have $T_{m}M_{X}={\rm Ker}\, D^2\psi_{|_{m}}$.
 \end{lem}

\noindent
{\bf Proof of Lemma \ref{lemH2}.}

The idea is to use a family of transformations which leave
$\overset{\circ}{M}_{X}$ and $n^2$ invariant (in a sense we define later), next to transport the equality $\text{Ker}\, D^2\psi_{|_{m_{0}}}= T_{m_{0}}M_{X}$  to any $m\in\overset{\circ}{M}_{X}$. 

\medskip
\mbox{}
\qquad
{\bf Family of transformations.} 
Let $m=(t,q,p,x,\xi,y,\eta)\in\overset{\circ}{M}_{X}$.
We write $m=(T_R,0,p,0,p,0,-p)$ for some $p\in\sqrt{2n^2(0)}\,\Sp^{d-1}$
 Thus, there exists $R_{p}\in\Ot(\R^d)$ such that $R_{p}(p)=p_{0}$. We define the map $\widetilde R_{m}:\R^{6d+1}\longrightarrow\R^{6d+1}$ by
 \begin{equation}
\label{rm}
\widetilde R_{m}(t,q,p,x,\xi,y,\eta)=
\left(t,R_{p}(q),R_{p}(y),R_p(x),R_{p}(\xi),R_p(y),R_{p}(\eta)\right).
 \end{equation}
 By construction we have $\widetilde{R}_{m}\,(m)=m_{0}$. 

\medskip
\mbox{}
\qquad
{\bf Action on the tangent place.} 
We have identified that the set $M_{X}$ satisfies
\begin{align*}
&
M_X
=
\underbrace{
\{(t,q,p,x,\xi,y,\eta), \text{ s.t. } \,
t=T_R, \, q=x=y=0, \, p=\xi=-\eta, \, p^2/2=n^2(0)\}
}_{
:=\widetilde{M_X}
}
\\
&
\qquad\qquad
\cap
\{p=(r,\theta_1,\ldots,\theta_{d-1})\, \text{ with } |\theta_1|\leq\theta_0\}.
\end{align*}
The set $\widetilde{M_X}$ is clearly invariant under the action of $\widetilde R_{m}$.
Therefore, by restricting the domain in the variable $\theta_1$, it is clear that whenever
$m\in \overset{\circ}{M}_{X}$, there exists a neighbourhood 
$U$ of $m$ in $\overset{\circ}{M}_{X}$ such that  $U_{0}:=\widetilde{R}_{m}U\subset\overset{\circ}{M}_{X}$. 
Since the application $\widetilde R_{m}$ is a linear map from $U$ to $U_{0}$ which satisfies $\widetilde R_{m}\, (m)=m_{0}$,
we deduce
$$
\widetilde{R}_{m}\ (T_{m}M_{X})=T_{m_{0}}M_{X}.
$$

\medskip
\mbox{}
\qquad
{\bf Action on the kernel.}  We now compute the set $\widetilde{R}_{m}(\text{Ker}(D^2\psi_{|_{m}}))$,
as follows
  \begin{align*}
 \widetilde{R}_{m}(\text{Ker}(D^2\psi_{|_{m}}))&=\{(T,R_{p}\,Q,R_{p}\,P,R_{p}\,X, R_p \, \Xi,R_{p}\,Y,R_{p}\,H),\ \text{s.t.}\ X=Y=Q=0,
\\
 &
\qquad
p.H=0,\ B_{T_{R}}(0,p) \, P+Tp=0,\ D_{T_{R}}(0,p) \, P=H\},
\\
&
=\{(T,Q,P,X,\Xi,Y,H),\ \text{s.t.}\ X=Y=Q=0,
\\
&
\qquad
p.R_{p}^{-1}H=0,\ B_{T_{R}}(0,p) \,R_{p}^{-1}P+Tp=0,\ D_{T_{R}}(0,p) \, R_{p}^{-1}P=R_{p}^{-1}H\}.
\\
 &
=\{(T,Q,P,X,\Xi,Y,H),\ \text{s.t.}\ X=Y=Q=0,
\\ 
 &
\qquad
p_{0}.H=0,\ R_{p}B_{T_{R}}(0,p) \, R_{p}^{-1}P+Tp_{0}=0,\ R_{p}D_{T_{R}}(0,p) \, R_{p}^{-1}P=H\}.
 \end{align*}
On the other hand, we claim that
\begin{equation}\label{eqstructBD}
 R_{p}B_{T_{R}}(0,p)R_{p}^{-1}=B_{T_{R}}(0,p_{0}),\qquad R_{p}D_{T_{R}}(0,p)R_{p}^{-1}=B_{T_{R}}(0,p_{0}).
 \end{equation}
Assuming the above identity is proved, we immediately deduce
$$
\widetilde{R}_{m}\, (\text{Ker}(D^2\psi_{|_{m}}))=\text{Ker} \, D^2\psi_{|_{m_{0}}}.
$$
We conclude by writing
$$
\widetilde{R}_{m}\, (\text{Ker}(D^2\psi_{|_{m}}))=\text{Ker}D^2\psi_{|_{m_{0}}}
=
T_{m_{0}}M_{X}
=
\widetilde{R}_{m}\ (T_{m}M_{X}).
$$
Thus, there only remains to prove (\ref{eqstructBD}).
By construction of the potential we clearly have
\begin{equation*}
R_p X(t,0,p)=X(t,0,p_0),
\text{ as well as }
n^2\left( R_p x) \right) =n^2\left( x\right),
 \end{equation*}
whenever $x/|x|$ lies in the angular sector $|\theta_1|\leq \theta_0$.
This  provides
$$
\frac{D^2n^2}{Dx^2}( X(t,0,p_{0}))
=
\frac{D^2n^2}{Dx^2}(R_{p}\, X(t,0,p))
=
R_{p}\, \frac{D^2n^2}{Dx^2}(X(t,0,p))\, R_{p}^{-1}.
$$
Therefore, using the differential equation (\ref{eqHJL2}) relating the time evolution of $B_t$ and $D_t$, we recover
the following system
\begin{equation*}
\label{eqHJL2}
\left\{
\begin{array}{ll}
\displaystyle
\vspace{0.2cm}
\frac{\partial}{\partial t}R_p B(t,0,p) R_p^{-1}=R_p D(t,0,p) R_p^{-1},& R_p B(0,0,p) R_p^{-1}=Id,\\
\displaystyle
 \frac{\partial}{\partial t} R_p D(t,0,p) R_p^{-1} =R_p \frac{D^2n^2}{Dx^2}(X(t,0,p)) B(t,0,p) R_p^{-1},
&
\\
\displaystyle
\qquad\qquad
= 
\frac{D^2n^2}{Dx^2}(R_p X(t,0,p)) R_p B(t,0,p) R_p^{-1}
&
\\
\displaystyle
\qquad\qquad
= 
\frac{D^2n^2}{Dx^2}( X(t,0,p_0)) R_p B(t,0,p) R_p^{-1}
&R_p D(0,0,p) R_p^{-1}=0.
\end{array}
\right.
\end{equation*}
Uniqueness of solutions to a differential system then gives
 \begin{equation*}
\forall t, \quad
R_{p}B_{t}(0,p)R_{p}^{-1}=B_t(0,p_0),
\quad
R_{p}D_{t}(0,p)R_{p}^{-1}=D_t(0,p_0).
 \end{equation*}
 Relation \eqref{eqstructBD} is proved.
 
\vspace{-0,3cm}
\hfill$\blacksquare$

\subsubsection{A useful byproduct of the proof of Proposition \ref{h2}}
\label{proofcorocoro}

\begin{lem}
\label{h22}
Let $n^2$ be the refraction index defined in \eqref{n2pot}.
Take any $m \in M_X$, written as $m=(T_R,0,p,0,p,0,-p)$ with
$p=\left(\sqrt{2 n^2(0)},\theta_1,\theta_2,\ldots,\theta_{d-1}\right)$
according to  Lemma \ref{lemstatset}.
Then,
$$
\psi(m) \text{ is constant on the set } |\theta_1|\leq \theta_0.
$$
\end{lem}

\medskip

\noindent
{\bf Proof of Lemma \ref{h22}.} 
Considering the actual value of $\psi(m)$, various terms need to be considered.
The term $\displaystyle \int_0^t (p_s^2/2+n^2(q_s)) \, ds$ is clearly constant whenever $|\theta_1|\leq \theta_0$.
The same statement holds for the factor $p_t \cdot q_t$.
The only non-obvious factor is $\Gamma_t q_t \cdot q_t$. As in the preceding proof we write
\begin{align*}
&
\Gamma_t(0,p) q_t(0,p) \cdot q_t(0,p)
=
\Gamma_t(0, p) q_t(0,R_p^{-1} p_0) \cdot
q_t(0,R_p^{-1} p_0)
\\
&
\qquad
=
R_p \Gamma_t(0,p) R_p^{-1} q_t(0, p_0) \cdot
q_t(0, p_0).
\end{align*}
There remains to write
\begin{align*}
&
R_p \Gamma_t(0,p) R_p^{-1}
=
R_p  \left( C_t(0,p) + i D_t(0,p) \right) \cdot \left( A_t(0,p) + i B_t(0, p) \right)^{-1}
R_p^{-1}
\\
&
\qquad
=
\left( R_p  C_t(0,p) R_p^{-1} + i R_p D_t(0,p) R_p^{-1} \right) \cdot
\left( R_p  A_t(0,p) R_p^{-1}  + i R_p  B_t(0, p) R_p^{-1} \right)^{-1}
\\
&
\qquad
=
\Gamma_t(0,p_0) 
\end{align*}
for we already know that $R_p  B_t(0, p) R_p^{-1} =B_t(0,p_0)$, 
$R_p  D_t(0, p) R_p^{-1} =D_t(0,p_0)$, and a similar proof establishes
$R_p  A_t(0, p) R_p^{-1} =B_t(0,p_0)$, 
$R_p  C_t(0, p) R_p^{-1} =D_t(0,p_0)$.

\vspace{-0,3cm}
\hfill$\blacksquare$

\subsection{The stationary phase argument:  Proof of item (iii) of our main Theorem}
\label{proofstat}

The main result of the present section is
\begin{prop}\label{lemcontributionimp}
Let $n^2$ be the potential constructed according to \ref{n2pot}.  Select a source $S\in\Sc(\R^d)$.
Then, the following holds.
\medskip

\noindent
(i)
If ${\rm supp}\left( \widehat{S}(\xi)\right) \cap \partial I_{\theta_0}=\emptyset$, we have
\begin{equation*}
\forall\,\phi\in\Sc(\R^d),\quad
\langle\widetilde{w_\eps}-L_\eps,\phi\rangle
=O_{T_1,\delta}(\sqrt\varepsilon),
\end{equation*} 
where $\langle L_{\varepsilon }, \phi\rangle$ is defined in (\ref{leps}) above (see also the Remark after Theorem \ref{propnonconvergence}),
and
$\partial I_{\theta_0}=\{\xi=(|\xi|,\theta_1,\ldots,\theta_{d-1}) \, \text{ such that } \, \theta_1=\pm \theta_0\}$
(see definition \ref{ith}).

\medskip

\noindent
(ii)
In the general case we have
\begin{equation*}
\forall\,\phi\in\Sc(\R^d),\quad
\langle\widetilde{w_\eps}-L_\eps,\phi\rangle
=o_{T_1,\delta}(\eps^0).
\end{equation*} 
\end{prop}

\medskip

\noindent
{\bf Proof of Proposition \ref{lemcontributionimp}.} 
Due to the fact that the stationary set $M_X$ in the to-be-developped stationary phase argument 
has a boundary at $\theta_1=\pm \theta_0$, the argument is  in two steps. This is the reason why the above
Proposition distinguishes between two cases.

\medskip

\noindent{\bf
$\bullet\bullet$ Proof of Proposition \ref{lemcontributionimp}-part (i)}

\medskip

Outside the stationary set $M_X$
associated with the complex phase $\psi$, the oscillatory integral~(\ref{oscinteg}) defining
$\langle\widetilde{w_\eps},\phi\rangle$
 is of order $O(\varepsilon^{\infty})$.
On the stationary set $M_X$ and near the support of $a_{N}$,  the stationary set $M_X$ is a submanifold {\em without boundary}, having codimension $k=6d+1-(d-1)=5d+2$. Indeed, thanks to the hypothesis on the support of $\widehat{S}$, we have $\text{supp}\, a_{N}\cap \partial M_{X}=\emptyset$.

\medskip

Let us now come to the explicit application of the stationary phase Theorem to the oscillatory integral~(\ref{oscinteg}).
Writing $p=(r,\theta_{1},\ldots,\theta_{d-1})$ in hyperspherical coordinates, we define the application:
 \begin{eqnarray*}
\gamma:\R^{6d+1}\cap \text{supp}\ a_N&\longrightarrow&\R^{5d+2}\times{\mathbb S}^{d-1}\\
(t,q,p,x,\xi,y,\eta)&\longmapsto& (\underbrace{t-T_{R},q,x,y,\xi-p,\eta+p,r-\sqrt{2n^2(0)}}_{=:\alpha},\underbrace{\theta_{1},\ldots,\theta_{d-1}}_{=:\theta})
\end{eqnarray*}
The map $\gamma$ is a $C^\infty$-diffeomorphism between $\text{supp}\ a_{N}$ and
$\gamma\left(\text{supp}\ a_{N}\right)$.
Furthermore, we have by construction
$$
(t,X)\in M_{X}\cap \text{supp}\ a_{N}\Longleftrightarrow\alpha=0.
$$
The new coordinates $(\alpha,\theta)$
 are adapted to the stationary set $M_{X}$ associated with $\psi$. Making the change of variables $(t,X)=\gamma^{-1}(\alpha,\theta)$ in the integral defining $\langle\widetilde{w_\eps},\phi\rangle$ we have
\begin{align}
\label{bientotfini}
&
\langle\widetilde{w_\eps},\phi\rangle
=
O_{\delta,T_{1}}(\varepsilon^N)
+
\\
&
\nonumber
\frac{1}{\varepsilon^{(5d+2)/2}}\int_{\gamma(supp\ a_{N})}
e^{\frac{i}{\varepsilon}\psi\circ\gamma^{-1}(\alpha,\theta)}
\left(\widehat{S}(.)\widehat{\phi}^*(.)P_{N}\left(.,.,.,\frac{.}{\sqrt\varepsilon}\right)\right)
\circ\gamma^{-1}(\alpha,\theta)\chi_{3}(\alpha,\theta)\, r^{d-1} \, d\alpha\, d\sigma(\theta),
\end{align}
where $d\sigma(\theta)$ denotes the standard euclidean surface measure on the unit sphere ${\mathbb S}^{d-1}$,
and $\chi_{3}$ is a truncation function on some compact set, a neighbourhood of $M_X$,
whose precise value is irrelevant. Here we have used the {\em non-stationary} phase Theorem
to reduce the original integral to an integral on a given compact set.

Since for all point $ m\in M_{X}\cap \text{supp}\, a_{N}$ we have $\text{Ker}(D^2\psi_{|_{m}})=T_{m}M_{X}$ (Lemma \ref{lemH2}), the function $D^2\psi$ is non-degenerate in the normal direction to $M_{X}$, which gives
\begin{equation*}
\text{det}\left(\frac{D^2\psi\circ\gamma^{-1}}{D\alpha^2}(0,\theta)\right)\neq0.
\end{equation*}
Furthermore, the projection of $\gamma(\text{supp}\ a_{N})$ onto the space variable $\theta$ is the angular sector
\begin{equation*}
\Pi_\theta I_{\theta_0}:=
\left\{(\theta_{1},\ldots,\theta_{d-1}),\quad \theta_{1}\in]-\theta_{0},\theta_{0}[\right\},
\end{equation*}
where $\Pi_\theta$ denotes the projection $(r,\theta_1,\ldots,\theta_{d-1}) \mapsto (\theta_1,\ldots,\theta_{d-1})$.
We can now apply the stationary phase Theorem in (\ref{bientotfini}).
Remembering that the codimension of the stationary set $M_{X}$ associated with $\psi$ is $5d+2$,
we obtain that for {\em any} integer $L$ there exists a sequence $(Q_{2\ell}(\partial))_{\ell\in\{0,\ldots,L\}}$
of operators of order $2\ell$ such that
\begin{align}
\label{eqJepsfin}
&
\nonumber
\langle\widetilde{w_\eps},\phi\rangle
=
\\
&
\nonumber
C_1 \, \int_{\Pi_{\theta}I_{\theta_{0}}}
\,
\frac{ \exp\left( i \frac{\pi}{4} {\rm sgn} \frac{D^2 \psi\circ \gamma^{-1}}{D \alpha^2}(0,\theta)\right)}
{\left| {\rm det} \frac{D^2 \psi\circ \gamma^{-1}}{D \alpha^2}(0,\theta)\right|}
\,
\exp\left(\frac{i}{\varepsilon}\psi\circ\gamma^{-1}(0,\theta)\right) \,
\\
&
\nonumber
\qquad\qquad\qquad\qquad\qquad
\left(\left(Q_{0}(.)\widehat S(.)\widehat{\phi^*}(.)P_{N}\left(.,.,.,\frac{.}{\sqrt{\varepsilon}}\right)\right)\circ\gamma^{-1}\chi_{3}\right)(0,\theta)
\, d\sigma(\theta)
\\
&
\nonumber
+
\int_{\Pi_{\theta}I_{\theta_{0}}}\exp\left(\frac{i}{\varepsilon}\psi\circ\gamma^{-1}(0,\theta)\right)\sum_{\ell=1}^L\varepsilon^\ell \, Q_{2\ell}(\partial)
\left(\left(\widehat S(.)\widehat{\phi^*}(.)P_{N}\left(.,.,.,\frac{.}{\sqrt{\varepsilon}}\right)\right)
\circ\gamma^{-1}\chi_{3}\right)(0,\theta) \,
d\sigma(\theta)
\\
&
+
O\left(\varepsilon^{L+1}\sup_{K\leq2L+d+3}\left\| \partial^K_{(\alpha,\theta)}
\left(\widehat S(.)\widehat{\phi^*}(.)P_{N}\left(.,.,.,\frac{.}{\sqrt{\varepsilon}}\right)\circ\gamma^{-1}\chi_{3}\right)\right\|_{L^\infty}\right)
+O_{\delta,T_{1}}(\varepsilon^N)
\\
&
\nonumber
:=I_\eps+II_\eps+III_\eps+O_{\delta,T_{1}}(\varepsilon^N),
\end{align}
with the value
$$
C_1=(2\pi)^{(5d+2)/2} \, (2 n^2(0))^{(d-1)/2}.
$$
The last line in
(\ref{eqJepsfin}) serves as a definition of the three terms $I_\eps$, $II_\eps$ and $III_\eps$,
and the $L^\infty$-norm in $III_\eps$ is evaluated on a compact set of values of $(\alpha,\theta)$,
whose precise  value is irrelevant.

We compute these three contributions.
Note that the retained value of the integer  $L$ remains to be determined at this stage.

\medskip

$\bullet$ \textbf{Contribution of  the remainder term $III_\eps$ in (\ref{eqJepsfin}).}

This term is best studied by coming back to the original variables $(t,X)$ instead of $(\alpha,\theta)$.
Expanding the $k$-th order derivatives involved in this term, we clearly have
\begin{align*}
&
III_\eps
=
O\left(\varepsilon^{L+1}\sup_{K\leq2L+d+3}\left\| \partial^K_{(t,X)}\left(\widehat S(.)\widehat{\phi^*}(.)P_{N}\left(.,.,.,\frac{.}{\sqrt{\varepsilon}}\right)\right)\right\|_{L^\infty}\right)
\\
&
\quad
=
O\left(\varepsilon^{L+1}\sup_{K\leq2L+d+3}\left\| \partial^K_{(t,X)} P_{N}\left(.,.,.,\frac{.}{\sqrt{\varepsilon}}\right)
\right\|_{L^\infty}\right).
\end{align*}
Hence, since
$$
P_{N}(t,q,p,x)=\pi^{-d/4}\text{det}(A(t,q,p)+iB(t,q,p))_{c}^{-1/2}\mathcal Q_{N}(t,q,p,x),
$$
we recover
\begin{align*}
&
III_\eps
=
O\left(\varepsilon^{L+1}\sup_{K\leq2L+d+3}\left\| \partial^K_{(t,q,p,y)} \left(
\mathcal Q_N(t,q,p,(y-q_t)/\sqrt{\eps})
\right)\right\|_{L^\infty}\right).
\end{align*}
Lastly, using \eqref{eqQN} we have
$$
\mathcal Q_{N}(t,q,p,x):=1+\sum_{(k,j)\in I_{N}}\varepsilon^{\frac{k}{2}-j}p_{k,j}(t,q,p,x),
$$
where $p_{k,j}$ has at most degree $k$ in $x$. We deduce
\begin{align*}
&
III_\eps
=
\sum_{(k,j)\in I_{N}}
O\left(\varepsilon^{\frac{k}{2}-j+L+1}
\sup_{K\leq2L+d+3}\left\| \partial^K_{(t,q,p,y)} \left(
p_{k,j}(t,q,p,(y-q_t)/\sqrt{\eps})
\right)\right\|\right)
\\
&
\quad
=
\sum_{(k,j)\in I_{N}}
O\left(\varepsilon^{\frac{k}{2}-j+L+1-\frac{k}{2}}\right)
=
O\left(\varepsilon^{L+1-(2N-1)}\right),
\end{align*}
where we have used that $j\leq 2N-1$ whenever $(k,j)\in I_N$ (see (\ref{eqQN})).
There remains to chose
$$
L=2N-1
$$
to recover
$$
III_\eps=O(\eps).
$$

$\bullet$ \textbf{Contribution of $II_\eps$ in (\ref{eqJepsfin}).}

This estimate is more delicate. Firstly, we have
\begin{align*}
&
II_\eps
=
\sum_{\ell=1}^L\varepsilon^\ell
O\left(
\left\|
 Q_{2\ell}(\partial)
\left(\left(\widehat S(.)\widehat{\phi^*}(.)P_{N}\left(.,.,.,\frac{.}{\sqrt{\varepsilon}}\right)\right)
\circ\gamma^{-1}\chi_{3}\right)(0,\theta)
\right\|_{L^\infty}
\right).
\end{align*}
Hence, going back to the $(t,X)$ variables again, and remembering that the relation
$(\alpha,\theta)=(0,\theta)$
implies $y=q_t=0$ and $t=T_R$, we recover the identity
\begin{align*}
&
II_\eps
=
\sum_{\ell=1}^L\varepsilon^\ell \, 
O\left(
\sup_{K\leq 2\ell}
\left\|
\partial^K_{(t,q,p,y)}\Big|_{y=q_t=0, t=T_R}
\left(
P_{N}\left(t,q,p,\frac{y-q_t}{\sqrt{\varepsilon}}\right)
\right)
\right\|_{L^\infty}
\right),
\end{align*}
where the $L^\infty$-norm is evaluated on some compact set of values of $p$.
Now, inserting the exact value of $P_N$, we may write
\begin{align*}
&
II_\eps
=
\sum_{\ell=1}^L\varepsilon^\ell \, 
O\left(
\sum_{(k,j)\in I_N}
\sup_{K\leq 2\ell}
\left\|
\partial^K_{(t,q,p,y)}\Big|_{y=q_t=0, t=T_R}
\left(
\eps^{\frac{k}{2}-j} \, p_{k,j}\left(t,q,p,\frac{y-q_t}{\sqrt{\varepsilon}}\right)
\right)
\right\|_{L^\infty}
\right)
\\
&
\quad
=
\sum_{\ell=1}^L \sum_{(k,j)\in I_N}
\varepsilon^\ell \, 
\eps^{\frac{k}{2}-j} \, 
O\left(
\sup_{K\leq 2\ell}
\left\|
\partial^K_{(t,q,p,y)}\Big|_{y=q_t=0, t=T_R}
\left(
p_{k,j}\left(t,q,p,\frac{y-q_t}{\sqrt{\varepsilon}}\right)
\right)
\right\|_{L^\infty}
\right)
\end{align*}
Hence, using the fact that each $p_{k,j}$ is a polynomial in its last argument, so that the above derivatives evaluated at
$y=q_t=0$ only leave the zero-th order term in the derived polynomial,
we recover
\begin{align*}
&
II_\eps
=
O\left(
\sum_{\ell=1}^L \sum_{(k,j)\in I_N}
\varepsilon^\ell \, 
\eps^{\frac{k}{2}-j} \, 
\sup_{K\leq 2\ell}
\eps^{-K/2}
\right)
\\
&
\quad
=
O\left(
\sum_{\ell=1}^L \sum_{(k,j)\in I_N}
\varepsilon^\ell \, 
\eps^{\frac{k}{2}-j} \, 
\eps^{-\ell}
\right)
=
O\left(
\sum_{\ell=1}^L \sum_{(k,j)\in I_N}
\eps^{\frac{k}{2}-j} \, \eps^{-\ell}
\right)
=
O\left(
\eps^{1/2}
\right),
\end{align*}
where we have used that $k-2j\geq 1$ whenever $(k,j)\in I_N$.
 
\medskip

$\bullet$ \textbf{Contribution of $I_\eps$ in (\ref{eqJepsfin}).}

The integral defining $I_\eps$ has the following more explicit value,
where $p=(\sqrt{2 n^2(0)},\theta_1,\ldots,\theta_{d-1})$, namely
\begin{multline*}
I_{\varepsilon}
=
C_1 \, \int_{\Pi_\theta I_{\theta_0}}
\frac{\displaystyle e^{i\frac \pi 4  \text{sgn}\left(\frac{D^2\psi\circ\gamma^{-1}}{D\alpha^2}(0,\theta)\right)}}{\text{det}\left(\frac{D^2\psi\circ\gamma^{-1}}{D\alpha^2}(0,\theta)\right) }
\,
\exp(\left(\frac{i}{\varepsilon}\psi(T_R,0,p,0,p,0,-p)\right)
\\
\text{det}(A(T_{R},0,p)+iB(T_{R},0,p))_{c}^{-1/2}\, \widehat S(p)\, \widehat{\phi}^*(-p)\, d\theta_{1}\ldots\,  d\theta_{d-1},
\end{multline*}
On top of that, we have
\begin{equation*}
\psi(T_R,0,p,0,p,0,-p)=\int_{0}^{T_{R}}\left(\frac{|p_{s}(0,p)|^2}{2}+n^2(q_{s}(0,p))\right) \, ds,
\end{equation*}
while the fact that $n^2$ is radial implies
that $\psi(T_R,0,p,0,p,0,-p)=\psi(T_R,0,p_0,0,p_0,0,-p_0)$ whenever $p\in I_{\theta_{0}}$.
For the same reason, we also have whenever $\theta\in\Pi_\theta I_{\theta_0}$ the relation
\begin{align*}
&
\frac{\displaystyle e^{i\frac \pi 4  \text{sgn}\left(\frac{D^2\psi\circ\gamma^{-1}}{D\alpha^2}(0,\theta)\right)}}{\text{det}\left(\frac{D^2\psi\circ\gamma^{-1}}{D\alpha^2}(0,\theta)\right) }
=
\frac{\displaystyle e^{i\frac \pi 4  \text{sgn}\left(\frac{D^2\psi\circ\gamma^{-1}}{D\alpha^2}(0,0)\right)}}{\text{det}\left(\frac{D^2\psi\circ\gamma^{-1}}{D\alpha^2}(0,0)\right) }
\end{align*}
together with the identity, valid whenever $p \in I_{\theta_0}$,
\begin{align*}
&
\text{det}(A(T_{R},0,p)+iB(T_{R},0,p))_{c}^{-1/2}
=
\text{det}(A(T_{R},0,p_0)+iB(T_{R},0,p_0))_{c}^{-1/2}.
\end{align*}
Eventually, we have obtained
\begin{equation}
\label{eqlim}
I_{\varepsilon}
=
C_{n^2,d} \,
e^{\left(\displaystyle\frac{i}{\varepsilon}\displaystyle\int_{0}^{T_{R}}\left(\frac{|p_{s}(0,p_0)|^2}{2}
+n^2(q_{s}(0,p_0))\right)ds\right)} \,
\int_{I_{\theta_{0}}}\widehat S(p)\widehat{\phi^*}(-p) \, d\sigma_{\theta_0}(p),
\end{equation}
with
\begin{equation*}
C_{T_{R},d}:=\frac{(2\pi)^{5d+2} e^{i\frac \pi 4 \text{sgn}\left(\frac{D^2\psi\circ \gamma^{-1}}{D\alpha^2}(0,0)\right)}}{\text{det}\left(\frac{D^2\psi\circ\gamma^{-1}}{D\alpha^2}(0,0)\right) }
\,
\text{det}(A(T_{R},0,p_{0})+iB(T_{R},0,p_{0}))_{c}^{-1/2}.
\end{equation*}
This ends the proof of Proposition \ref{lemcontributionimp}-part (i).

\medskip

\noindent{\bf
$\bullet\bullet$ Proof of Proposition \ref{lemcontributionimp}-part (ii)}

\medskip

In that case, the argument is essentially the same (a stationary phase argument in the variable~$\alpha$),
up to a convenient use of the
dominated convergence Theorem (to deal with the variable $\theta_1$, and more specifically with the boundary
$\theta_1=\pm \theta_0$).

Namely, we first write, as in the proof of part (i) of the Proposition,
\begin{align}
\label{bientotfini2}
&
\langle\widetilde{w_\eps},\phi\rangle
=
O_{\delta,T_{1}}(\varepsilon^N)
+
\\
&
\nonumber
\frac{1}{\varepsilon^{(5d+2)/2}}\int_{\gamma(supp\ a_{N})}
e^{\frac{i}{\varepsilon}\psi\circ\gamma^{-1}(\alpha,\theta)}
\left(\widehat{S}(.)\widehat{\phi}^*(.)P_{N}\left(.,.,.,\frac{.}{\sqrt\varepsilon}\right)\right)
\circ\gamma^{-1}(\alpha,\theta)\chi_{3}(\alpha,\theta)\, r^{d-1} \, d\alpha\, d\sigma(\theta),
\end{align}
where $\chi_{3}$ is a truncation function on some compact set, a neighbourhood of $M_X$,
whose precise value is irrelevant. Here we have used the non-stationary phase Theorem
to reduce the original integral to an integral on a given compact set.
The key point now lies in writing, 
\begin{align}
\label{bientotfini3}
&
\langle\widetilde{w_\eps},\phi\rangle
=
O_{\delta,T_{1}}(\varepsilon^N)
+
\\
&
\nonumber
\int d\sigma(\theta) \, 
\underbrace{
\left(
\frac{1}{\varepsilon^{(5d+2)/2}}
\int d\alpha \, 
e^{\frac{i}{\varepsilon}\psi\circ\gamma^{-1}(\alpha,\theta)}
\left(\widehat{S}(.)\widehat{\phi}^*(.)P_{N}\left(.,.,.,\frac{.}{\sqrt\varepsilon}\right)\right)
\circ\gamma^{-1}(\alpha,\theta)\chi_{3}(\alpha,\theta)\, r^{d-1}\right)
}_{=:J_\eps(\theta)}.
\end{align}
With this formulation in mind, our next objective is to prove that whenever $\eta>0$ is a small parameter we have
\begin{align}
\label{bientotfini4}
&
\int_{|\theta_1 \pm \theta_0|\leq \eta} d\sigma(\theta) \, 
\left| J_\eps(\theta) \right| \leq C \, \eta,
\end{align}
for some $C>0$ independent of $\eps$ and $\eta$.
It is clear indeed that the upper-bound (\ref{bientotfini4}), in conjunction with part (i) of the Proposition,
provides a complete proof of Proposition \ref{lemcontributionimp}-part~(ii).

Let us now concentraate on the case $|\theta_1-\theta_0|\leq \eta$
(the proof in the case $|\theta_1+\theta_0|\leq \eta$ is the same).

In order to prove (\ref{bientotfini4}), 
we fix a value  $(\theta_2^0,\ldots,\theta_{d-1}^0)$ and we prove that,
given $(\theta_2^0,\ldots,\theta_{d-1}^0)$, there is an $\eta>0$, and a $C>0$ independent of $\eps$, such that
\begin{align}
\label{bientotfini5}
&
\forall \theta \, \text{ such that } \, |\theta-(\theta_0.\theta_2^0,\ldots,\theta_{d-1}^0)|\leq \eta, \quad
\text{ we have } \,
\left| J_\eps(\theta) \right| \leq C.
\end{align}
Covering the whole set $\{\theta \in {\mathbb S}^{d-1} ; |\theta_1-\theta_0|\leq \eta\}$ by finitely many sets of the 
form $\{ |\theta-(\theta_0.\theta_2^0,\ldots,\theta_{d-1}^0)|\leq \eta\}$ clearly provides the desired relation
(\ref{bientotfini4}) once (\ref{bientotfini5}) is proved.

Now, relation (\ref{bientotfini5})
results from an application of the stationary phase Theorem, with {\em complex phase} and  {\em with parameter}.
Here $\alpha$ 
is the variable used for the stationary phase itself, while $\theta$ is the parameter, and $\psi\circ\gamma^{-1}$ is the
complex phase. We introduce the short-hand notation
$\theta^0=(\theta_0,(\theta')^{0})=(\theta_0,\theta_2^0,\ldots,\theta_{d-1}^0)$ for convenience.
It has already been established\footnote{{\em Stricto sensu}, these relations have only be proved 
when $|\theta_1|<\theta_0$, and we here extend the result to the case $\theta_1=\theta_0$.
This is allowed due to the invariance of the phase on the parameter $\theta$ whenever $|\theta_1|\leq \theta_0$ --
Lemma \ref{h22}.}
that
\begin{align*}
&
{\rm Im}\left(\psi\circ\gamma^{-1}\right)(\alpha,\theta) \geq 0, \quad \forall (\alpha,\theta),
\\
&
{\rm Im} \left(\psi\circ\gamma^{-1}\right)(\alpha=0,\theta=\theta^0)=0,
\\
&
\nabla_\alpha\left(\psi\circ\gamma^{-1}\right)(\alpha=0,\theta=\theta^0)=0,
\\
&
{\rm det}
\left(
\frac{
D^2 \psi\circ\gamma^{-1}
}
{D \alpha^2}
\right)(\alpha=0,\theta=\theta^0)
\neq
0
\end{align*}
Therefore, the stationary phase theorem with parameter ensures
that close to $\theta=\theta^0$ there is an expansion of the form
\begin{align*}
J_\eps(\theta)
=
e^{i \phi(\theta)/\eps}
\,
\left(
\sum_{\ell=0}^L
\eps^\ell \, 
\left(
Q_{2\ell}(\partial_\alpha) 
\left(
\widehat{S}(.) \,
\widehat{\phi^*}(.) \,
P_N\left( ., ., \frac{.}{\sqrt{\eps}}\right)
\circ\gamma^{-1} \,
\chi_3(.)
\right)
\right)^0(\theta)
\right)
+
R(\eps,L,\theta),
\end{align*}
for some smooth functions $\phi$ and $R(\eps,L,\theta)$, where the $Q_{2\ell}$'s are differential operators of order $2\ell$
in the variable $\alpha$, and, for any function $u(\alpha,\theta)$, the notation $u^0(\theta)$ refers to any smooth
function $u^0(\theta)$ that belongs to the same residue class than the original function $u(\alpha,\theta)$ modulo
the ideal generated by $\nabla_\alpha \psi \circ \gamma^{-1}(\alpha,\theta)$
(see H\"ormander  \cite{MR1996773}, sect. 7.7, for the details). With this notation, we actually have
$\phi=\left(\psi\circ\gamma^{-1}\right)^0$.
Besides, the remainder term $R$ satisfies
as the term $III_\eps$ in the previous step an estimate of the form
\begin{align*}
\left|
R(\eps,L,\theta)
\right|
\leq
C_L \, \varepsilon^{L+1} \,
\left(
\sup_{K\leq 2 (L+1)}
\left\|
\partial^K_\alpha
\left(
\left(
\widehat S(.) \, 
\widehat{\phi^*}(.) \, 
P_{N}\left(.,.,.,\frac{.}{\sqrt{\varepsilon}}\right) \,
\right)
\circ\gamma^{-1} \,
\chi_3(.)
\right)
\right\|_{L^\infty}
\right),
\end{align*}
for some constant $C_L>0$ independent of $\eps$, and provided  $\theta$ is close to $\theta^0$ (independently of $\eps$).
These two ingredients immediately provide, using the same estimates as we did for the terms $III_\eps$ and $II_\eps$
above, the upper-bound, valid for $\theta$ close to $\theta^0$,
\begin{align*}
\left|
J_\eps(\theta)
\right|
\leq
C \, 
\left(
\sum_{\ell=0}^L
\eps^\ell \, 
\left(
Q_{2\ell}(\partial_\alpha) 
\left(
\widehat{S}(.) \,
\widehat{\phi^*}(.) \,
P_N\left( ., ., \frac{.}{\sqrt{\eps}}\right)
\circ\gamma^{-1} \,
\chi_3(.)
\right)
\right)^0(\theta)
\right)
+
R(\eps,L,\theta),
\end{align*}
Gathering powers of $\eps$ as in the previous part of the proof, provides the upper bound
\begin{align*}
\left|
J_\eps(\theta)
\right|
\leq
C,
\end{align*}
where $C$ does not depend on $\eps$ and $\theta$ is close to $\theta^0$, independently of $\eps$.
Point (\ref{bientotfini5}) is proved.

We immediately deduce that (\ref{bientotfini4}) holds, and the proof of
Proposition \ref{lemcontributionimp} -- part (ii) is complete.

\hfill$\blacksquare$\\


\subsection{Conclusion}


Gathering the intermediate result in Proposition \ref{intermres}, together with 
Proposition \ref{lemcontributionimp}, gives item (iii)
of  Theorem \ref{propnonconvergence},
by 
 conveniently choosing the parameters $\delta$, $\theta$, $T_{0}$ and $T_{1}$.

\bibliographystyle{plain}\bibliography{bibliothese}

\begin{thebibliography}{10}

\bibitem{MR1924691}
Jean-David Benamou, Fran{\c{c}}ois Castella, Theodoros Katsaounis, and Benoit
  Perthame.
\newblock High frequency limit of the {H}elmholtz equations.
\newblock {\em Rev. Mat. Iberoamericana}, 18(1):187--209, 2002.

\bibitem{MR2582438}
Jean-Fran{\c{c}}ois Bony.
\newblock Mesures limites pour l'\'equation de {H}elmholtz dans le cas non
  captif.
\newblock {\em Ann. Fac. Sci. Toulouse Math. (6)}, 18(3):459--493, 2009.

\bibitem{MR1900556}
F.~Castella, B.~Perthame, and O.~Runborg.
\newblock High frequency limit of the {H}elmholtz equation. {II}. {S}ource on a
  general smooth manifold.
\newblock {\em Comm. Partial Differential Equations}, 27(3-4):607--651, 2002.

\bibitem{MR2139886}
Fran{\c{c}}ois Castella.
\newblock The radiation condition at infinity for the high-frequency
  {H}elmholtz equation with source term: a wave-packet approach.
\newblock {\em J. Funct. Anal.}, 223(1):204--257, 2005.

\bibitem{MR2289541}
Fran{\c{c}}ois Castella and Thierry Jecko.
\newblock Besov estimates in the high-frequency {H}elmholtz equation, for a
  non-trapping and {$C\sp 2$} potential.
\newblock {\em J. Differential Equations}, 228(2):440--485, 2006.

\bibitem{CJK}
Fran{\c{c}}ois Castella, Thierry Jecko, and Andr\'eas Knauf.
\newblock Semiclassical resolvent estimates for schrödinger operators with
  coulomb singularities.
\newblock {\em Ann. H. Poincaré}, 9(4):775--815, 2008.

\bibitem{MR1461126}
M.~Combescure and D.~Robert.
\newblock Semiclassical spreading of quantum wave packets and applications near
  unstable fixed points of the classical flow.
\newblock {\em Asymptot. Anal.}, 14(4):377--404, 1997.

\bibitem{MR1735654}
Mouez Dimassi and Johannes Sj{\"o}strand.
\newblock {\em Spectral asymptotics in the semi-classical limit}, volume 268 of
  {\em London Mathematical Society Lecture Note Series}.
\newblock Cambridge University Press, Cambridge, 1999.

\bibitem{MR2237164}
Elise Fouassier.
\newblock High frequency analysis of {H}elmholtz equations: case of two point
  sources.
\newblock {\em SIAM J. Math. Anal.}, 38(2):617--636 (electronic), 2006.

\bibitem{MR2296804}
Elise Fouassier.
\newblock High frequency limit of {H}elmholtz equations: refraction by sharp
  interfaces.
\newblock {\em J. Math. Pures Appl. (9)}, 87(2):144--192, 2007.

\bibitem{MR724029}
B.~Helffer and D.~Robert.
\newblock Calcul fonctionnel par la transformation de {M}ellin et op\'erateurs
  admissibles.
\newblock {\em J. Funct. Anal.}, 53(3):246--268, 1983.

\bibitem{MR1996773}
Lars H{\"o}rmander.
\newblock {\em The analysis of linear partial differential operators. {I}}.
\newblock Classics in Mathematics. Springer-Verlag, Berlin, 2003.
\newblock Distribution theory and Fourier analysis, Reprint of the second
  (1990) edition [Springer, Berlin; MR1065993 (91m:35001a)].

\bibitem{MR1872698}
Andr{{\'e}} Martinez.
\newblock {\em An introduction to semiclassical and microlocal analysis}.
\newblock Universitext. Springer-Verlag, New York, 2002.

\bibitem{MR1695559}
Benoit Perthame and Luis Vega.
\newblock Morrey-{C}ampanato estimates for {H}elmholtz equations.
\newblock {\em J. Funct. Anal.}, 164(2):340--355, 1999.

\bibitem{MR897108}
Didier Robert.
\newblock {\em Autour de l'approximation semi-classique}, volume~68 of {\em
  Progress in Mathematics}.
\newblock Birkh\"auser Boston Inc., Boston, MA, 1987.

\bibitem{FixMe}
Julien Royer.
\newblock {Limiting absorption principle for the dissipative Helmholtz
  equation.}
\newblock {\em Commun. Partial Differ. Equations}, 35(8):1458--1489, 2010.

\bibitem{MR2718664}
Julien Royer.
\newblock Semiclassical measure for the solution of the dissipative {H}elmholtz
  equation.
\newblock {\em J. Differential Equations}, 249(11):2703--2756, 2010.

\bibitem{MR927007}
Xue~Ping Wang.
\newblock Time-decay of scattering solutions and resolvent estimates for
  semiclassical {S}chr\"odinger operators.
\newblock {\em J. Differential Equations}, 71(2):348--395, 1988.

\bibitem{MR2184186}
Xue~Ping Wang and Ping Zhang.
\newblock High-frequency limit of the {H}elmholtz equation with variable
  refraction index.
\newblock {\em J. Funct. Anal.}, 230(1):116--168, 2006.

\end{thebibliography}
\end{document}